\documentclass[lettersize,journal]{IEEEtran}

\usepackage{amsmath,amsfonts,mathtools}
   \usepackage{empheq}

\usepackage{algorithmic}
\usepackage{algorithm}
\usepackage{array}
\usepackage[caption=false,font=normalsize,labelfont=sf,textfont=sf]{subfig}
\usepackage{textcomp}
\usepackage{stfloats}
\usepackage{url}
\usepackage{verbatim}
\usepackage{graphicx}
\usepackage{cite}
\usepackage[dvipsnames]{xcolor}
\usepackage{circuitikz}

\usepackage{amsthm}
\newtheoremstyle{break}
{\topsep}{\topsep}%
{\itshape}{}%
{\bfseries}{}%
{\newline}{}%
\theoremstyle{break}

\providecommand{\tx}[1]{\text{\upshape{#1}}}	
\providecommand{\ddt}[1]{\dfrac{\tx{d}{#1}}{\tx{d}t}}
\usepackage[printonlyused]{acronym}
\acrodef{hac}[HAC]{hybrid angle control}
\acrodef{coi}[COI]{center-of-inertia}
\acrodef{ib}[IB]{infinite bus}
\acrodef{sg}[SG]{synchronous generators}
\acrodef{wrt}[w.r.t.]{with respect to}
\acrodef{agas}[AGAS]{almost global asymptotic stability}
\acrodef{lhs}[LHS]{left-hand side}  
\acrodef{rhs}[RHS]{right-hand side}  
\acrodef{rocof}[RoCoF]{rate of change of frequency}  
\acrodef{gfm}[GFM]{grid-forming}  
\acrodef{gfl}[GFL]{grid-following}  
\acrodef{chil}[C-HiL]{controller-hardware-in-the-loop}
\acrodef{sm}[SM]{synchronous machine}
\acrodef{pcc}[PCC]{point of common coupling}
\begin{document}
\definecolor{col1}{RGB}{214,45,32} 
\definecolor{col2}{RGB}{0,87,231} 
\definecolor{col3}{RGB}{255,167,0} 
\definecolor{col4}{RGB}{0,135,68} 
\definecolor{col5}{RGB}{179, 179, 179} 

\title{Grid-Forming Hybrid Angle Control: \\Behavior, Stability, Variants and Verification}

\author{Ali Tayyebi, Denis Vettoretti, Adolfo Anta, and Florian D\"{o}rfler 
	\thanks{A. Tayyebi (the corresponding author) is with the Hitachi Energy Research (HER), 72226 V\"{a}sterås, Sweden, and Automatic Control Laboratory, ETH Z\"{u}rich, 8092 Z\"{u}rich, Switzerland,  e-mail: ali.tayyebi@hitachienergy.com, D. Vettoretti and A. Anta are with the Austrian Institute of Technology (AIT), 1210 Vienna, Austria and F. D\"{o}rfler is with the Automatic Control Laboratory, ETH Z\"{u}rich, 8092 Z\"{u}rich, Switzerland. This work was funded in parts by the {HER}, AIT, and ETH funds.}}

\maketitle

\begin{abstract}
This work explores the stability, behavior, variants, and a \acf{chil} verification of the recently proposed \acf{gfm} \acf{hac}. We revisit the foundation of \ac{gfm} \ac{hac}, and highlight its behavioral properties in relation to the conventional \ac{sm}. Next, we introduce the required complementary controls to be combined with the \ac{hac} to realize a \ac{gfm} behavior. The characterization of the analytical operating point and nonlinear energy-based stability analysis of a grid-connected converter under the \ac{hac} is presented. Further, we consider various output filter configurations and derive an approximation for the original control proposal. Moreover, we provide details on the integration of \ac{gfm} \ac{hac} into a complex converter control architecture and introduce several variants of the standard \ac{hac}. Finally, the performance of \ac{gfm} \ac{hac} is verified by several test scenarios in a \ac{chil} setup to test its behavior against real-world effect such as noise and delays.
\end{abstract}


\section{Introduction}
The global shift toward the massive integration of energy generation from renewable source accompanied by the supply chain concerns associated with conventional energy generation has raised significant interest in converter-based systems. Thus, power converters are perceived as the vital corner stones of the modern power system and are expected to replace the well-established \ac{sm} technology. However, a robust and reliable control of power converters in a converter-dominated power system is to some extent an open question. The emerging \acf{gfm} control synthesis in contrast to the classic \acf{gfl} converter control concept is envisioned to address the stability challenges in a converter-dominated power grid \cite{hatziargyriou2020definition,lin2020research,crivellaro2020beyond,milano2018foundations,TGAKD20,Rocabert,markovic2021understanding}. On the other hand, it is worth mentioning that the power system operators are actively designing test procedures and grid code requirements for the \ac{gfm} converters, as well, e.g.,  \cite{NGESO2023,AEMO2022}. 

The broadly recognized droop control serves as a powerful baseline \ac{gfm} control candidate that mimics the behavior of a \ac{sm} governor for the power converters \cite{CDA93,rowe2012arctan,yu2020comparative}. As the natural extension of the droop control, the virtual synchronous machine concept is proposed  that emulates the \ac{sm} dynamics (up to different degrees of accuracy) \cite{ZW11,chen2021enhanced}. On the other hand, the matching-type \ac{gfm} controllers are proposed that synthesize the converter control based on the structural dynamic similarities with the \ac{sm} \cite{arghir2019electronic,cvetkovic_modeling_2015}. Along a different design direction, the nonlinear oscillators dynamics are recently exploited for a \ac{gfm} control design \cite{aghdam2022virtual,seo2019dispatchable,awal2022double}. Finally, the combination of aforementioned techniques has resulted in several hybrid control architectures \cite{tayyebi2022grid,gao2020grid,gross2022dual,chen2022generalized}. 

The \ac{gfm} \ac{hac} relies on a combination of the dc matching control and ac synchronization term that resembles the droop control and/or Kuramoto oscillator dynamics \cite{simpson2013synchronization}. The theoretic control design and system-level simulation-based performance investigation of the \ac{hac} are previously explored \cite{tayyebi2022hybrid,tayyebi2020almost,tayyebi2022system,tayyebi2022grid}. Previous works highlight 1) the strong stability properties of the \ac{hac} under mild parametric conditions, 2) system-level frequency stability enhancement, 3) stabilizing behavior in complex hybrid ac/dc power grids, and finally, 4) robustness \ac{wrt} the nonlinear phenomena such as current limitation and grid split. 

In this paper, we highlight further details on the behavioral properties of the \ac{gfm} \ac{hac}, provide a closed-form characterization of the converter operating point, formulate an intuitive energy-based stability analysis, provide an approximate form of the \ac{hac}, and disclose several variants of the standard \ac{hac}. Last but not least,  the performance of \ac{gfm} \ac{hac} is verified by several test scenarios in a \ac{chil} setup that utilizes an OPAL-RT simulator and external control cards.    

The remainder of this paper is structured as it follows. Section \ref{TSG-sec: stiff grid} describes the dynamic modeling of a grid-connected converter, introduces the \ac{hac}, and discusses its behavioral properties. Section \ref{TSG-sec:stability} presents the closed-loop analysis. Section \ref{TSG-sec: weak grid} provides details on the weak grid connection, \ac{hac} approximation, and its variants. Section \ref{TSG-sec:verification} provides \ac{chil} performance verification, and Section \ref{TSG-sec:conclusion} concludes the paper.
\section{Converter connected to a stiff grid}\label{TSG-sec: stiff grid}
In this section, we present the dynamic modeling of a grid-connected converter, revisit the \ac{gfm} \ac{hac} strategy, and discuss its behavioral properties.
\subsection{Physical converter system dynamics}\label{TSG-subsec: dynamic model}
Let us consider a two-level dc-ac power converter model that is interfaced to a stiff grid (i.e., with constant frequency and voltage magnitude) through a resistive-inductive element \cite{arghir2019electronic}; see Figure \ref{TSG-fig:1}. The open-loop dc voltage and ac current dynamics of such system are described by
\begin{subequations}\label{TSG-eq:1}
\begin{align}
C_\tx{dc}\ddt{v_\tx{dc}}&=i_\tx{dc}-G_\tx{dc}v_\tx{dc}-i_\tx{s},
\label{TSG-eq:1a}
\\
L\ddt{i_\tx{abc}}&=v_\tx{s,abc}-Ri_\tx{abc}-v_\tx{g,abc},
\label{TSG-eq:1b}
\end{align}
\end{subequations}
where $C_\tx{dc}$ denotes the dc-link capacitance, $v_\tx{dc}$ denotes the dc-link voltage, $i_\tx{dc}$ denotes the current flowing out of the dc energy source, $G_\tx{dc}$ denotes the dc conductance that models the dc-side losses, and $i_\tx{s}$ denotes the dc-side switched current. Further, $L$ denotes the ac-side inductance that models the converter output filter, $i_\tx{abc}$ denotes the converter output current flowing into the grid, $v_\tx{s,abc}$ denotes the ac-side switched voltage, $R$ denotes the equivalent series resistance of the filter inductance, and finally $v_\tx{g,abc}$ denotes the balanced voltage of the stiff ac grid at nominal frequency $\omega_\tx{0}$ and magnitude $v_\tx{0}$. 
\subsection{Energy source model and control}\label{TSG-subsec: energy source model}
\begin{figure}[t!]
	\centering	
	{\includegraphics[width=0.8\columnwidth]{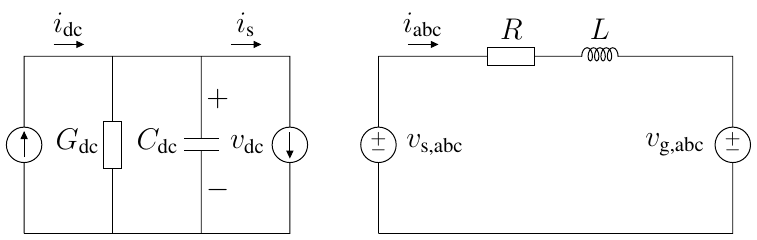}}
	\caption{The circuit diagram associated with the open-loop dynamics of the grid-connected converter model as in \eqref{TSG-eq:1}.}\label{TSG-fig:1}
\end{figure}

The dc current source in Figure \ref{TSG-fig:1} that models the primary dc energy source e.g., a battery, can be controlled in several ways. One can close the loop by considering a proportional controller to increase the dc voltage damping \cite{TGAKD20}. Further, it is possible to also include integral and derivative terms to enhance the dc voltage reference tracking and dynamic performance. Thus, the $i_\tx{dc}$ in \eqref{TSG-eq:1a} takes generic form
\begin{equation}\label{TSG-eq:2}
i_\tx{dc}=-\kappa_\tx{p} \left(v_\tx{dc}-v_\tx{dc,r}\right)-\kappa_\tx{i}\int_{0}^{t} \left(v_\tx{dc}-v_\tx{dc,r}\right) \tx{d}\tau-\kappa_\tx{d}\ddt{v_\tx{dc}},
\end{equation}  
where $\kappa_\tx{p}$, $\kappa_\tx{i}$, and $\kappa_\tx{d}$ denote the parameters of the proportional-integral-derivative control. Note that one can implement the derivative term in \eqref{TSG-eq:2} by measuring the dc-link capacitance current (since they are equivalent up to a constant factor).  Previous work investigated the contribution of control \eqref{TSG-eq:2} to the frequency damping and inertial response under the \ac{gfm} matching control \cite{arghir2019electronic,arghir2018grid}. Finally, if the energy source is not controllable, e.g., battery-integrated system without a dc-dc converter stage, one can fix $i_\tx{dc}$ to a constant reference.
\subsection{Power-preserving averaged dc-ac converter model} \label{TSG-subsec: averaged model}
The dc-ac converter \eqref{TSG-eq:1} is represented by the switched current and voltage pair $\left(i_\tx{s},v_\tx{s,abc}\right)$. Let us introduce the balanced three-phase converter modulation signal 
\begin{equation}\label{TSG-eq:3}
m_\tx{abc}=\mu\left[\cos\theta,\cos\left(\theta+\dfrac{2\pi}{3}\right),\cos\left(\theta-\dfrac{2\pi}{3}\right)\right]^\top,
\end{equation}
where $\mu$ and $\theta$ respectively denote the modulation signal magnitude and phase angle. Next, the lossless power-preserving averaged model of the two-level converter is given by \cite{yazdani10voltage}
\begin{subequations}\label{TSG-eq:4} 
\begin{align}
i_\tx{s}&=m_\tx{abc}^\top i_\tx{abc},
\\
v_\tx{s,abc}&=v_\tx{dc}m_\tx{abc}.
\end{align}
\end{subequations}
In the sequel, we show how $\mu$ and $\theta$ are selected. 
Aside from the nonlinear power-preserving model \eqref{TSG-eq:4}, other approaches can be considered. 
\begin{itemize}
\item One can assume decoupled dc and ac dynamics that is usually verified by considering a sufficiently fast dc voltage control \cite{CDA93,ZW11,chen2021enhanced}. In this viewpoint, the converter is seen as an ideal controllable voltage source.
\item Another trend is to model the internal dynamics of the dc-ac converters by oscillator dynamics. This approach is often adopted when studying the stability of interconnected converter-based systems; see \cite{aghdam2022virtual} for a review and \cite{seo2019dispatchable,yu2020comparative} for experimental investigations.
\item Recent works highlight the application of hybrid systems theory in modeling the converter dynamics. These works consider a blend of discontinuous and continuous signals in converter dynamical description, therefore, do not distinguish between the switching and continuous averaged converter models; e.g., see \cite{colon2022stability,albea2017hybrid}. 
\end{itemize}  
\subsection{Grid-forming hybrid angle control strategy}
In this subsection, we briefly revisit the design of grid-forming \ac{hac} \cite{tayyebi2020almost,tayyebi2022grid,tayyebi2022hybrid,tayyebi2022system}. Let us begin by defining the converter relative angle \ac{wrt} the grid model in Figure \ref{TSG-fig:1}. The modulation angle $\theta$ in \eqref{TSG-eq:3} enters $v_\tx{s,abc}$  in \eqref{TSG-eq:4} that subsequently appears in \eqref{TSG-eq:1b} (as the voltage behind output filter). Let us define $\omega$ as the converter angular frequency that is given by the time-derivative of $\theta$. Similarly, let $\theta_\tx{g}$ and $\omega_\tx{g}$ respectively denote the phase angle and the angular frequency grid voltage $v_\tx{g,abc}$ in \eqref{TSG-eq:1b} and Figure \ref{TSG-fig:1}. Note that $\omega_\tx{g}=\omega_0$, since we consider a stiff grid. Thus, the converter-grid relative angle and its derivative are given by
\begin{subequations}
\begin{align}
{\delta}&=\theta-\theta_\tx{g},
\label{TSG-eq:5a}
\\
\ddt{\delta}&=\ddt{\theta}-\ddt{\theta_\tx{g}}=\omega-\omega_\tx{g}=\omega-\omega_0.
\label{TSG-eq:5b}
\end{align}
\end{subequations}   
The \ac{gfm} \ac{hac} that defines the converter frequency (hence the modulation angle) takes the form
\begin{equation}\label{TSG-eq:6}
\boxed{\omega=\omega_0+\underbracket{\kappa_\tx{dc}\left(v_\tx{dc}-v_\tx{dc,r}\right)}_\tx{dc matching term}-\underbracket{\kappa_\tx{ac}\sin\left(\dfrac{\delta-\delta_\tx{r}}{2}\right)}_\tx{ac synchronization term},}
\end{equation}
where $v_\tx{dc,r}$ and $\delta_\tx{r}$ respectively denote the dc voltage and relative angle references. It is worth mentioning that the dc part of \eqref{TSG-eq:6} is similar to the matching control \cite{arghir2019electronic,cvetkovic_modeling_2015,CGD17}. On the other hand, the ac part of the \ac{hac} realizes the frequency synchronization via nonlinear angle damping assignment. The prior works \cite{tayyebi2022grid,tayyebi2022system,TGAKD20} provide detailed discussions on the properties of \ac{gfm} controls that depend on the ac and/or dc quantities. Nonetheless, in a nutshell, incorporating the dc feedback in the frequency dynamics tends to enhance the robustness and including the ac feedback enhances the dynamic performance; see \cite{gao2020grid,tayyebi2022system,samanta2022fast,samanta2021stability} for theoretic and numerical investigations.  The \ac{hac} \eqref{TSG-eq:6}, while defining the converter frequency, behaves as a synchronization mechanism. To further elaborate, if the converter dc voltage is sufficiently regulated, i.e., $v_\tx{dc} \approx v_\tx{dc,r}$, then the converter-grid relative angle dynamics \eqref{TSG-eq:5b} reduces to  
\begin{equation}\label{TSG-eq:7}
\ddt{\delta}\approx-\kappa_\tx{ac}\sin\left(\dfrac{\delta-\delta_\tx{r}}{2}\right).
\end{equation}
This means  if $\delta>\delta_\tx{r}~\Rightarrow~\tx{d}\delta/\tx{d}t<0~\Rightarrow~\delta\downarrow$ and similarly, if $\delta<\delta_\tx{r}~\Rightarrow~\tx{d}\delta/\tx{d}t>0~\Rightarrow~\delta\uparrow$.
The HAC potentially replaces the synchronization mechanism (e.g., phase-locked loop, virtual synchronous machine, active power control sub-systems) in converter control architectures; see Figure \ref{TSG-fig:5}. For instance, HAC is a synchronizing control candidate for the
\begin{itemize}
\item high voltage direct current (HVDC) converters in embedded, inter-connector, multi-terminal, and offshore wind farm integration setups,
\item flexible ac transmission system (FACTS) devices,
\item low-voltage photovoltaic (PV) and battery systems,
\item and, utility-scale battery energy storage system (BESS).
\end{itemize}
\subsection{DC voltage and AC power flow regulation}
The \ac{hac} regulates the dc voltage and ac power flow through frequency synchronization. To further elaborate, let us consider two separate cases. 

\subsubsection{Pure dc feedback control} assume $\kappa_\tx{dc}\neq0$ and $\kappa_\tx{ac}=0$ that reduces \eqref{TSG-eq:6} to
\begin{equation}
	\Delta\omega=\omega-\omega_0=\kappa_\tx{dc}\left(v_\tx{dc}-v_\tx{dc,r}\right)=\kappa_\tx{dc}\Delta v_\tx{dc}.
\end{equation}
This controller combination is the reduction of \ac{hac} to the matching control \cite{arghir2019electronic}. It is established that under the matching control, i.e., when the converter frequency is defined proportional to the dc voltage, the converter dynamics are structurally similar to that of the SM. Therefore, the converter exhibits self-synchronizing behavior of the \ac{sm} \cite{cvetkovic_modeling_2015} which means
\begin{equation*}
	\omega\to\omega_0 \Rightarrow \Delta\omega\to 0 \Rightarrow  \Delta v_\tx{dc}\to 0 \Rightarrow v_\tx{dc}\to v_\tx{dc,r}.
\end{equation*} 
Thus, the frequency synchronization implies dc voltage regulation, that is achieved by modifying the ac power. This control mode is particularly interesting in weak dc-link applications.

\subsubsection{Pure ac feedback control} consider the gain combination $\kappa_\tx{dc}=0$ and $\kappa_\tx{ac}\neq0$.  Let us approximate the ac term in \eqref{TSG-eq:6} with the ac power flow deviation, i.e., assume that $\Delta\delta=\delta-\delta_\tx{r}$ is proportional to $\Delta p=p-p_\tx{r}$ up to a constant factor $\kappa_{\delta-p}$. Then, \eqref{TSG-eq:6} reduces to 
\begin{equation}\label{TSG-eq:droop control}
	\Delta\omega\approx-\kappa_\tx{ac}\sin\left(\dfrac{\kappa_{\delta-p}\Delta p}{2}\right)\approx-\left(\dfrac{\kappa_\tx{ac}\kappa_{\delta-p}}{2}\right)\Delta p,
\end{equation}
assuming that $\Delta p$ is sufficiently small. This variant represents the power-frequency droop control embedded in \ac{hac}, thus, 
\begin{equation*}
	\omega\to\omega_0 \Rightarrow \Delta\omega\to 0 \Rightarrow \Delta p\to 0\Rightarrow p\to p_\tx{r}.
\end{equation*}       	
In this case, frequency synchronization implies ac power flow regulation, that is achieved by the power injection/absorption of the dc-link. This control mode is particularly interesting in stiff dc-link applications. Finally, the hybrid configuration under appropriate tuning provides seamless transition between the aforementioned modes \cite{tayyebi2022grid,gross2022dual}.

\subsection{Behavioral interpretations and connections to the SM}
It is possible to interpret the structure of \ac{hac} in relation to \ac{sm} control and behavior. Firstly, the influence of a governor on the \ac{sm} behavior is perceived as modifying the turbine output mechanical power $p_\tx{m}$ according to the mechanical frequency $\omega_\tx{m}$ deviation from its reference $\omega_\tx{m,r}$. In other words,
\begin{equation}\label{TSG-eq:governor}
	p_\tx{m}=p_\tx{m,r}-\kappa_{\omega-p}\left(\omega_\tx{m}-\omega_\tx{m,r}\right),
\end{equation}
where $p_\tx{m,r}$ and $\kappa_{\omega-p}$ respectively denote the turbine reference power and governor control gain. Observe that if $\omega_\tx{m}\uparrow\downarrow~\Rightarrow~ p_\tx{m}\downarrow\uparrow$ to accordingly modify the energy input into the \ac{sm} such that the frequency is stabilized. One can alternatively rewrite \eqref{TSG-eq:governor} as the so-called droop control, i.e.,
\begin{equation}\label{TSG-eq: inverse governor}
	\omega_\tx{m}=\omega_\tx{m,r}-\dfrac{1}{\kappa_{\omega-p}}\left(p_\tx{m}-p_\tx{m,r}\right).
\end{equation}
Now, under the small power-angle assumption, i.e., $\Delta\delta\propto\Delta p$ one can interpret the ac part of the \ac{hac} \eqref{TSG-eq:6} as droop control \eqref{TSG-eq:droop control} which takes the same form as \eqref{TSG-eq: inverse governor}. Therefore, the ac term in \eqref{TSG-eq:6} mimics the stabilizing influence of the turbine governor.
 
Next, let us revisit the modeling of \ac{sm} inertial response \cite{ulbig2014impact}. We assume the that the mechanical power $p_\tx{m}$ is flowing into the \ac{sm} and electrical power $p_\tx{e}$ is flowing out of its ac terminal. These two quantities are linked through the time-derivative of kinetic energy $E_\tx{k}$ stored in the \ac{sm} rotor, i.e.,
 \begin{equation}
 	\ddt{E_\tx{k}}=p_\tx{m}-p_\tx{e}\quad\tx{where}\quad E_\tx{k}=\dfrac{1}{2}J\omega^2
 \end{equation}
 and $J$ denotes the rotor moment of inertia. A salient feature of the \ac{sm} is that if there is an imbalance between its mechanical and electrical powers, e.g., due to load variation, the rotating mass acts as an energy buffer and provides/absorbs the excess power to restore the power balance. The resulting influence is the \ac{sm} frequency variation, i.e.,
 \begin{equation}
 	\text{if}~\ddt{E_\tx{k}}=J\omega\ddt{\omega}>0~(\tx{or}~<0)~\Rightarrow~\omega\uparrow(\downarrow).
 \end{equation}
 The dc-ac power converters, by design, incorporate a similar mechanism. To further elaborate, let $p_\tx{dc}$ denote the power that is flown into the converter dc-link and $p_\tx{ac}$ is the power that is flown out of the converter ac terminal. These quantities are linked together through the potential energy $E_\tx{p}$ that is stored in the converter dc-link, i.e., 
 \begin{equation}
 	\ddt{E_\tx{p}}=p_\tx{dc}-p_\tx{ac}\quad\tx{where}\quad E_\tx{p}=\dfrac{1}{2}C_\tx{dc}v_\tx{dc}^2.
 \end{equation}
 Similarly, the power imbalance between the converter dc and ac ports is compensated by the dc-link energy variation, i.e.,
 \begin{equation}
 	\text{if}~\ddt{E_\tx{p}}=C_\tx{dc}v_\tx{dc}\ddt{v_\tx{dc}}>0~(\tx{or}~<0)~\Rightarrow~v_\tx{dc}\uparrow(\downarrow).
 \end{equation}
 From this perspective, the dc term in \ac{hac} \eqref{TSG-eq:6} that relates the converter frequency to the dc voltage (i.e., $\omega\propto v_\tx{dc}$), resembles the inertial response of the \ac{sm} and links the converter frequency to the available physical stored energy in the dc-link capacitance.
\subsection{AC voltage control}
The \ac{gfm} \ac{hac} is primarily designed as an active power-frequency controller \cite{tayyebi2022grid}. Thus, one has to consider complementary ac voltage control. Similar to other grid-forming controls \cite{TGAKD20}, there are different control candidates.
\begin{itemize}
\item One can implement a proportional-integral (PI) (or simply a proportional) ac voltage control that processes the \ac{pcc} voltage error and provides a reference magnitude for  converter modulation in \eqref{TSG-eq:3} \cite{tayyebi2022system}.
\item Another alternative is to define the modulation signal magnitude based on a reactive power and voltage droop control \cite{markovic2021understanding,Rocabert}. In this approach, the converter modulation magnitude is modified if the reactive power deviates form its reference. Thus, the modulation magnitude modification indirectly controls the \ac{pcc} voltage.
\item The most straightforward, although less robust, approach is to define the reference magnitude for the converter modulation signal according to the desired references for the dc and ac voltages \cite{tayyebi2022grid}.     
\end{itemize}

\subsection{Control implementation and filtering requirement}\label{TSG-subsec:implementation}
The previous work \cite{tayyebi2022system}, establish that \ac{hac} \eqref{TSG-eq:6} can be exactly constructed based on the dc voltage measurement, internal converter modulation angle, and the grid voltage measurement in Figure \ref{TSG-fig:1}. To recapitulate, one should firstly expand the ac term in \eqref{TSG-eq:6}, i.e.,
\begin{equation}
\sin\left(\dfrac{\delta-\delta_\tx{r}}{2}\right)=\sin\dfrac{\delta}{2}\cos\dfrac{\delta_\tx{r}}{2} - \cos\dfrac{\delta}{2}\sin\dfrac{\delta_\tx{r}}{2}.
\end{equation}
Then, the terms depending on $\delta_\tx{r}$ can be computed according to the prescribed power and voltage set-points \cite{tayyebi2022grid}. Next, the terms depending on $\delta=\theta-\theta_\tx{g}$ are constructed based on the sines an cosines of $\theta$ and $\theta_\tx{g}$ that can be respectively obtained from the converter modulation signal $m_\tx{abc}$ and the grid voltage $v_\tx{g,abc}$. Note that it is standard practice to low-pass filter the dc voltage feedback in \eqref{TSG-eq:6} and the grid voltage measurement to remove the potential dc ripple and ac noise, respectively.

\section{Closed-loop stability analysis}\label{TSG-sec:stability}
In what follows, we select a combination of the controls described in the previous section, construct the closed-loop dynamics, and investigate the overall system stability.
\subsection{Closed-loop system formulation}
Let us begin by transforming the three-phase dynamics \eqref{TSG-eq:1} to the stationary $\alpha\beta$-coordinates by using the standard Clarke transformation \cite{yazdani10voltage} that results in 
\begin{subequations}\label{TSG-eq:alphabeta-model}
	\begin{align}
		C_\tx{dc}\ddt{v_\tx{dc}}&=i_\tx{dc}-G_\tx{dc}v_\tx{dc}-i_{\tx{s}},
		\label{TSG-eq:alphabeta-model1}
		\\
		L\ddt{i_{\alpha\beta}}&=v_{\tx{s},\alpha\beta}-Ri_{\alpha\beta}-v_{\tx{g},\alpha\beta}.
		\label{TSG-eq:alphabeta-model2}
	\end{align}
\end{subequations}
Next, we select the PI dc voltage control\footnote{The previous works on \ac{hac} \cite{tayyebi2022grid,tayyebi2020almost,tayyebi2022hybrid,tayyebi2022system} do not include the integral term in their dc voltage controls, therefore, the forthcoming closed-loop system analysis (although conceptually similar) differs from the prior investigations. In particular, the PI dc voltage control consideration omits the previously required assumption to prove the existence and derive a closed-form expression of the closed-loop stationary operating points.} (for an enhanced dynamic performance and robustness) from \eqref{TSG-eq:2}, i.e.,
\begin{subequations}\label{TSG-eq:proportional dc voltage control}
\begin{align*}
\ddt{\zeta}
&=
v_\tx{dc}-v_\tx{dc,r},
\\
i_\tx{dc}
&=
-\kappa_\tx{p} \left(v_\tx{dc}-v_\tx{dc,r}\right)-\kappa_\tx{i}\zeta,
\end{align*}
\end{subequations}
where  $\zeta$ denotes the integrator state,
the \ac{hac} \eqref{TSG-eq:6}, and the feedforward ac voltage control, i.e.,
\begin{equation}\label{TSG-eq:feedforward AC voltage control}
\mu=\dfrac{v_\tx{r}}{v_\tx{dc,r}}.
\end{equation}  
Therefore, all three control inputs, i.e., the dc energy source current, modulation magnitude, and angle are well-defined (the latter is obtained by integrating the converter frequency defined by \eqref{TSG-eq:6}). Next, we consider rotating dq-coordinates \cite{yazdani10voltage} that are aligned with the grid angle $\theta_\tx{g}$, thus, rotating with the grid frequency $\omega_\tx{g}$. The closed-loop dynamics in rotating dq-coordinates is represented by      
\begin{subequations}\label{TSG-eq:CL system}
\begin{align}
\ddt{\delta} 
&=
\kappa_\tx{dc}\left(v_\tx{dc}-v_\tx{dc,r}\right)-\kappa_\tx{ac}\sin\left(\dfrac{\delta-\delta_\tx{r}}{2}\right),\label{TSG-eq:CL system1}\\
\ddt{\zeta}
&=
v_\tx{dc}-v_\tx{dc,r},\label{TSG-eq:CL system2}\\
C_\tx{dc}\ddt{v_\tx{dc}}
&=
-\kappa_\tx{p} \left(v_\tx{dc}-v_\tx{dc,r}\right)-\kappa_\tx{i}\zeta-G_\tx{dc}v_\tx{dc}\nonumber\\
&~~~-{\mu \left(i_\tx{d}\cos\delta+i_\tx{q}\sin\delta\right)},\label{TSG-eq:CL system3}\\
L\ddt{i_\tx{d}}
&=
{\mu v_\tx{dc}\cos\delta} - R i_\tx{d}-L\omega_0i_\tx{q}-v_\tx{g,d},\label{TSG-eq:CL system4}\\
L\ddt{i_\tx{q}}
&=
{\mu v_\tx{dc}\sin\delta} - R i_\tx{q}+L\omega_0i_\tx{d}.\label{TSG-eq:CL system5}
\end{align}
\end{subequations}
We remark that $v_\tx{g,d}=v_0$ (i.e., the nominal voltage magnitude of the stiff grid) and $v_\tx{g,q}=0$ since the d-axis is aligned with $\theta_\tx{g}$. Moreover, it is important to emphasize that the closed-loop system is nonlinear due to \ac{hac} in \eqref{TSG-eq:CL system1} and modulated trigonometric terms in \eqref{TSG-eq:CL system3} and \eqref{TSG-eq:CL system4}.
\subsection{Analytical derivation of equilibria}
In order to evaluate the stationary operating points (denoted by star superscript) of the closed-loop system, we begin by setting the \ac{rhs} of the \eqref{TSG-eq:CL system} to zero, i.e.,
\begin{subequations}\label{TSG-eq:SS system}
	\begin{align}
		\kappa_\tx{dc}\left(v_\tx{dc}^\star-v_\tx{dc,r}\right)-\kappa_\tx{ac}\sin\left(\dfrac{\delta^\star-\delta_\tx{r}}{2}\right)&=0,\label{TSG-eq:SS system1}\\
		v_\tx{dc}^\star-v_\tx{dc,r}&=0,\label{TSG-eq:SS system2}\\
		-\kappa_\tx{p} \left(v_\tx{dc}^\star-v_\tx{dc,r}\right)-\kappa_\tx{i}\zeta^\star-G_\tx{dc}v_\tx{dc}^\star&\nonumber\\
		-{\mu \left(i_\tx{d}^\star\cos\delta^\star+i_\tx{q}^\star\sin\delta^\star\right)}&=0,\label{TSG-eq:SS system3}\\
		{\mu v_\tx{dc}^\star\cos\delta^\star} - R i_\tx{d}^\star-L\omega_0i_\tx{q}^\star-v_\tx{g,d}&=0,\label{TSG-eq:SS system4}\\
		{\mu v_\tx{dc}^\star\sin\delta^\star} - R i_\tx{q}^\star+L\omega_0i_\tx{d}^\star&=0.\label{TSG-eq:SS system5}
	\end{align}
\end{subequations}
Hence, one can solve \eqref{TSG-eq:SS system1} and \eqref{TSG-eq:SS system2} to evaluate $\delta^\star$ and $v_\tx{dc}^\star$. Next, it is possible to solve \eqref{TSG-eq:SS system4} and \eqref{TSG-eq:SS system5} that result in the closed-form expressions for $i_\tx{d}^\star$ and $i_\tx{q}^\star$ that are the functions of $\delta^\star$ and $v_\tx{dc}^\star$. Finally, one can solve \eqref{TSG-eq:SS system3} to find $\zeta^\star$. Thus, (letting $k\in\{1,2\}$) the steady-state system of equations \eqref{TSG-eq:SS system} yields the following operating points,
\begin{subequations}\label{TSG-eq:equilibria}
\begin{align}
\delta^\star&=\delta_\tx{r}+2\pi k,\label{TSG-eq:equilibria1}
\\
v_\tx{dc}^\star&=v_\tx{dc,r},\label{TSG-eq:equilibria2}
\\
\zeta^\star&=\dfrac{-G_\tx{dc} v_\tx{dc}^\star-{\mu \left(i_\tx{d}^\star\cos\delta^\star-i_\tx{q}^\star\sin\delta^\star\right)}}{\kappa_\tx{i}},\label{TSG-eq:equilibria3}
\\
i_\tx{d}^\star&=\dfrac{\mu v_\tx{dc}^\star\left(R\cos\delta^\star-L\omega_0\sin\delta^\star\right)-Rv_\tx{g,d}}{R^2+\left(L\omega_0\right)^2},\label{TSG-eq:equilibria4}
\\
i_\tx{q}^\star&=\dfrac{\mu v_\tx{dc}^\star\left(L\omega_0\cos\delta^\star+R\sin\delta^\star\right)-L\omega_0v_\tx{g,d}}{R^2+\left(L\omega_0\right)^2}.\label{TSG-eq:equilibria5}
\end{align}	
\end{subequations}
We remark that due to the periodicity of \eqref{TSG-eq:equilibria3}-\eqref{TSG-eq:equilibria5} \ac{wrt} $\delta$, the steady-state quantities $\zeta^\star$, $i_\tx{d}^\star$, and $i_\tx{q}^\star$ are identical for either $\delta_\tx{r}$ or $\delta_\tx{r}+2\pi$ \cite{tayyebi2022grid}. In the next subsection, we investigate the stability of the operating point in \eqref{TSG-eq:equilibria} that is characterized by $\delta^\star=\delta_\tx{r}$, i.e.,
\begin{equation}\label{TSG-eq:desired equilibrium}
\boxed{x^\star=\left(\delta_\tx{r},\zeta^\star, v_\tx{dc}^\star, i_\tx{d}^\star, i_\tx{q}^\star\right).}
\end{equation}
In the sequel, we restrict our focus to a local state space region around the stationary point \eqref{TSG-eq:desired equilibrium} that excludes the other angle equilibrium in \eqref{TSG-eq:equilibria1}. The reader is referred to \cite{tayyebi2022grid} for a global (i.e., large-signal) stability analysis\footnote{We remark that the forthcoming analysis can be extended to to provide large-signal, i.e., global stability guarantees, as well. However, in this paper, for the sake of brevity, a local analysis is provided.}.
\subsection{Nonlinear energy-based stability analysis}
The system \eqref{TSG-eq:CL system} is characterized by the nonlinearities due to the \ac{hac} in \eqref{TSG-eq:CL system1}, and the modulated current and voltage terms in \eqref{TSG-eq:CL system3}-\eqref{TSG-eq:CL system5}. One potential approach to analyze the stability of the operating point \eqref{TSG-eq:desired equilibrium}, is to linearize the \eqref{TSG-eq:CL system} and investigate the eigenvalues of resulting linear system. However, due to the particular structure of the Jacobian associated with \eqref{TSG-eq:CL system}, it is not straightforward to derive the analytical closed-form expressions for the eigenvalues evaluated at \eqref{TSG-eq:desired equilibrium}. 

A more comprehensive nonlinear analysis approach is to associate a so-called \emph{energy function} with the closed-loop dynamics \eqref{TSG-eq:CL system} and study the behavior of this function \ac{wrt} the evolution of states in \eqref{TSG-eq:CL system}. The energy function behavior (under certain conditions) reveals the stability properties of the system. More precisely, let us define an energy function      
\begin{equation}\label{TSG-eq:LF}
	V(\hat{x})=c_1 \left(1-\cos\dfrac{\hat{\delta}}{2}\right)+c_2\hat{\zeta}^2+c_3 \hat{v}_\tx{dc}^2+c_4 \hat{i}_\tx{d}^2+c_5 \hat{i}_\tx{q}^2,
\end{equation}
\begin{figure*}[t!]
	\begin{align}\label{TSG-eq:LF derivative}
		\ddt{V(\hat{x})}&=\left(\dfrac{c_1}{2}\sin\dfrac{\hat{\delta}}{2}\right)\ddt{\hat{\delta}}+2\left(\left(c_2\hat{\zeta}\right)\ddt{\hat{\zeta}}+\left(c_3\hat{v}_\tx{dc}\right)\ddt{\hat{v}_\tx{dc}}+\left(c_4\hat{i}_\tx{d}\right)\ddt{\hat{i}_\tx{d}}+\left(c_5\hat{i}_\tx{q}\right)\ddt{\hat{i}_\tx{q}}\right),\nonumber
		\\
		&=\left(\dfrac{c_1}{2}\sin\dfrac{\delta-\delta^\star}{2}\right)\ddt{\delta}+2\left(c_2\left(\zeta-\zeta^\star\right)\ddt{\zeta}+c_3\left(v_\tx{dc}-v_\tx{dc}^\star\right)\ddt{v_\tx{dc}}+c_4\left(i_\tx{d}-i_\tx{d}^\star\right)\ddt{i_\tx{d}}+c_5\left(i_\tx{q}-i_\tx{q}^\star\right)\ddt{i_\tx{q}}\right).
	\end{align}
	\hrule
\end{figure*}
\noindent where $\hat{x}=x-x^\star$, $\hat{\delta}=\delta-\delta^\star$, $\hat{\zeta}=\zeta-\zeta^\star$, $\hat{v}_\tx{dc}=v_\tx{dc}-v_\tx{dc}^\star$, $\hat{i}_\tx{d}=i_\tx{d}-i_\tx{d}^\star$, and $\hat{i}_\tx{q}=i_\tx{q}-i_\tx{q}^\star$ and all the coefficients in \eqref{TSG-eq:LF} are positive constants. Let us consider the derivative of \eqref{TSG-eq:LF} \ac{wrt} time as in \eqref{TSG-eq:LF derivative}.  The state-dependent function \eqref{TSG-eq:LF} can be perceived as a measure of the distance (i.e., error) between the states in \eqref{TSG-eq:CL system} and the equilibrium point \eqref{TSG-eq:desired equilibrium}. 

\begin{figure}[t!]
	\centering	
	\includegraphics[trim=7mm 8mm 4mm 5mm,clip,width=0.8\columnwidth]{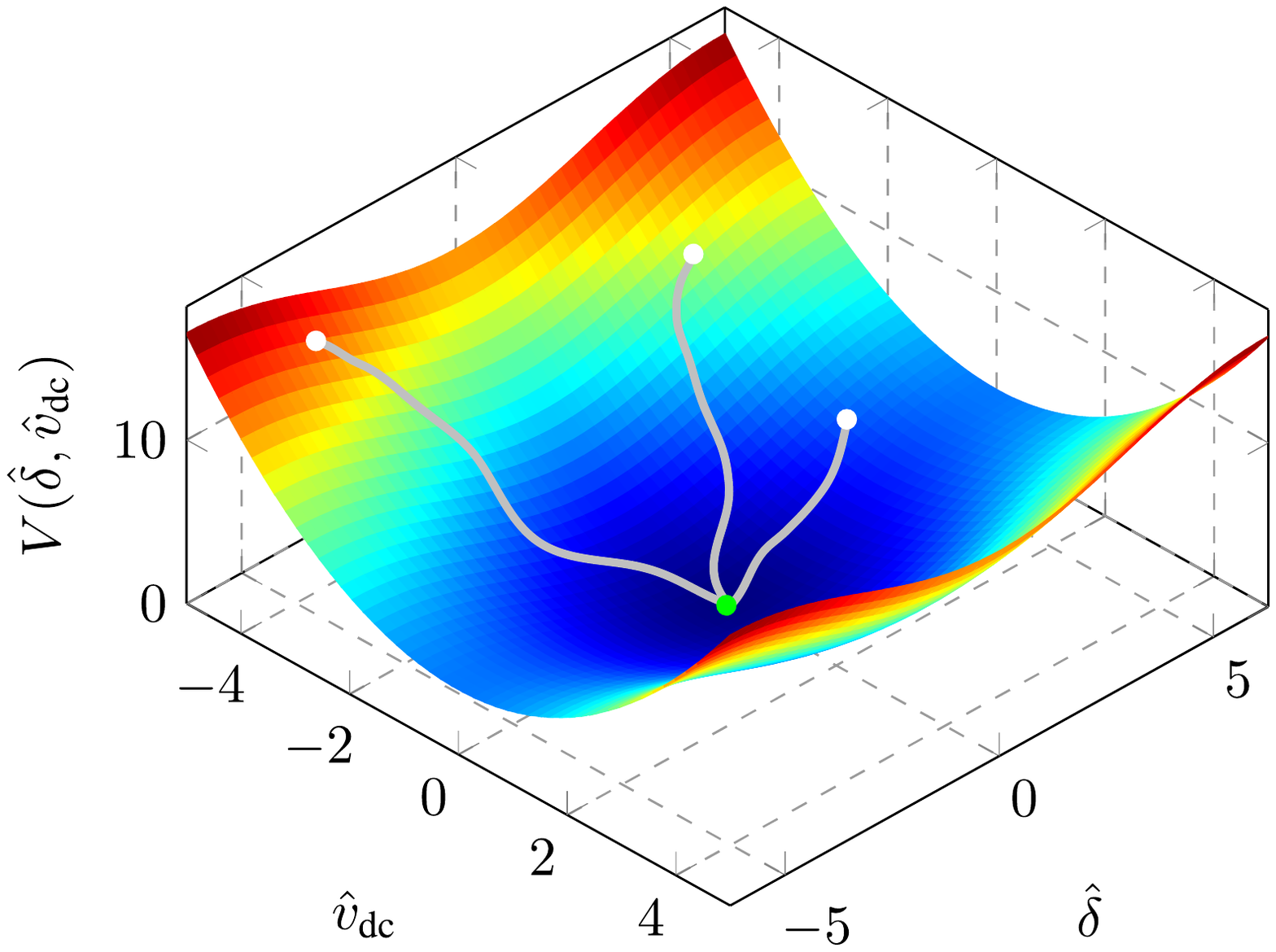}
	\caption{Conceptual illustration of the energy-based stability analysis; here it is assumed that the energy function is only function of the relative angle and dc voltage, i.e., $V\big(\hat{\delta},\hat{v}_\tx{dc}\big)=2\big(1-\cos\hat{\delta}/2\big)+(1/2)\hat{v}_\tx{dc}^2$. }\label{TSG-fig:2}	
\end{figure}

Given that $V(0)=0$ and $V(\hat{x})>0$ for $x\neq x^\star$, we can conclude that $V(\hat{x})\to 0$ (thus, $\hat{x}\to 0$ and $x\to x^\star$) \textit{if} $\tx{d}V(\hat{x})/\tx{d}t<0$ for all $x\neq x^\star$ (i.e., if $V(\hat{x})$ is strictly decreasing). Hence, the convergence of $V(\hat{x})$ to zero implies the stability of \eqref{TSG-eq:CL system}. In order to demonstrate the stability of system \eqref{TSG-eq:CL system}, we seek for a parametric condition that results in $\tx{d}V(\hat{x})/\tx{d}t<0$ for all $x\neq x^\star$. Following the procedure in \cite[Theorem 2]{tayyebi2022grid}, we set the coefficients in \eqref{TSG-eq:LF} as
\begin{equation*}
c_1=\dfrac{4}{\kappa_\tx{dc}}, c_2=\dfrac{\kappa_\tx{i}}{2}, c_3=\dfrac{C_\tx{dc}}{2},~\tx{and}~c_4=c_5=\dfrac{L}{2}.
\end{equation*} 
Subsequently, lengthy albeit straightforward  computation as in \cite[Theorem 2]{tayyebi2022grid}\cite{Andrea} shows that if
\begin{equation}\label{TSG-eq: condition}
\boxed{\rho=\dfrac{\kappa_\tx{ac}}{\kappa_\tx{dc}}>\rho_\tx{critical},}
\end{equation}
where 
\begin{equation*}
\rho_\tx{critical}=\dfrac{1}{G_\tx{dc}+\kappa_\tx{p}}+\dfrac{\mu^2\left(i_\tx{d}^{\star2}+i_\tx{q}^{\star2}\right)}{G_\tx{dc}+\kappa_\tx{p}}+\dfrac{\mu^2 v_\tx{dc}^{\star2}}{R},
\end{equation*}
then $\tx{d}V(\hat{x})/\tx{d}t<0$. We remark that the condition \eqref{TSG-eq: condition} is met by choosing a sufficiently large ac synchronization gain in \eqref{TSG-eq:6}. Note that the implication of the stability condition \eqref{TSG-eq: condition} is the fact that the closed-loop stability is guaranteed solely by an appropriate choice of the \ac{hac} ac and dc gains. Further, $\rho_\tx{critical}$ can be reduced by increasing the proportional gain of the dc voltage control, thus, allowing for a less aggressive tuning of the \ac{hac}; see \cite{tayyebi2022grid} for details on the stability condition \eqref{TSG-eq: condition}. Last, Figure \ref{TSG-fig:2} provides a conceptual example for the presented energy-based stability analysis where a simplified form of \eqref{TSG-eq:LF} is employed. Figure \ref{TSG-fig:2} shows that how different initial states converge to the desired equilibrium point if the energy decay, i.e., $\tx{d}V(\hat{x})/\tx{d}t<0$ is guaranteed.
\section{Converter connected to a weak grid}\label{TSG-sec: weak grid}
In this section, we consider more complex model configurations in contrast to the model presented in Figure \ref{TSG-fig:1}. Further, an approximate variant of \ac{hac} is presented. Next, we show how \ac{hac} can be combined with classic cascaded current and voltage controls. Finally, we present several \ac{hac} variants. 
\subsection{Grid impedance consideration}
\begin{figure}[b!]
	\includegraphics[width=\columnwidth]{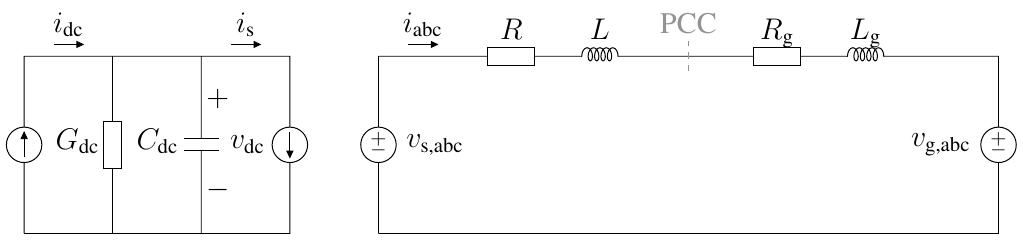}
\caption{The circuit diagram of the converter model connected to a weak grid model in abc-coordinates system.}\label{TSG-fig:3}
\end{figure}

A weak grid connection is considered by including an equivalent grid impedance that is represented by a resistive-inductive element as shown in Figure \ref{TSG-fig:3}. The different ratios of $L_\tx{g}$ and $R_\tx{g}$ represent connection to the low, medium, and high voltage grids \cite{Rocabert}. Since the serial connected filter and grid equivalent impedances can be merged together, the closed-loop dynamics associated with the model in Figure \ref{TSG-fig:3} takes the same form as in \eqref{TSG-eq:CL system}. The main implication of a weak grid connection is the fact that $v_\tx{g,abc}$ (in Figure \ref{TSG-fig:1}) is not available for the control implementation; see subsection \ref{TSG-subsec:implementation} and \cite{tayyebi2022grid}. Next, we show how leveraging certain assumptions allows to counteract this limitation by deriving an approximation for the HAC based on the ac active power flow.
\subsection{Power-based control approximation}
Consider the model in Figure \ref{TSG-fig:3} and let us merge the filter and grid equivalent impedance into a unified resistive-inductive element that reduces the ac sub-circuit to a classic coupled voltage sources configuration as in \cite[Figure 7]{Rocabert}. Subsequently, under dominantly inductive grid and small power angle assumptions \cite{Rocabert,rowe2012arctan}, the relative angle between $v_\tx{s,abc}$ and $v_\tx{g,abc}$ is linearly approximated by the active power flows, i.e.,
\begin{equation}\label{TSG-eq:approximation 1}
\delta \approx \sin(\delta) \approx \left(\dfrac{\left(L+{L}_\tx{g}\right)\omega_0}{|v_\tx{s,abc}||v_\tx{g,abc}|}\right) p,
\end{equation} 
where $p$ denotes the power injected by the converter. Further, assuming regulated ac voltages, i.e., constant $|v_\tx{s,abc}|$ and $|v_\tx{g,abc}|$, \eqref{TSG-eq:approximation 1} is simplified to
\begin{equation}\label{TSG-eq:approximation 2}
	\delta \approx \alpha  p\quad\tx{where}\quad \alpha = \dfrac{\left(L+{L}_\tx{g}\right)\omega_0}{|v_\tx{s,abc}||v_\tx{g,abc}|}\quad\tx{is constant}.
\end{equation} 
A similar approximation as in \eqref{TSG-eq:approximation 2} relates $\delta_\tx{r}$ to the power reference $p_\tx{r}$. Thus, the \ac{hac} in \eqref{TSG-eq:6} is approximated by
\begin{equation}\label{TSG-eq: approximation 3}
\boxed{\omega\approx\omega_0+\kappa_\tx{dc}\left(v_\tx{dc}-v_\tx{dc,r}\right)-\bar{\kappa}_\tx{ac}\left(p-p_\tx{r}\right),}
\end{equation} 
where $\bar{\kappa}_\tx{ac}=\alpha \kappa_\tx{ac}/2$. The approximate HAC \eqref{TSG-eq: approximation 3}, can be re-written in a trade-off form as
\begin{equation}\label{TSG-eq: approximation 4}
\Delta \omega \approx \kappa_\tx{dc} \Delta v_\tx{dc} -  \bar{\kappa}_\tx{ac} \Delta p,
\end{equation} 
where the converter frequency deviation from the nominal frequency is proportional to the dc voltage and ac power deviations from their respective references. We remark that if the assumptions behind \eqref{TSG-eq: approximation 3} hold, the local (i.e., small-signal) stability properties of the original and approximate HAC forms are identical. Finally, note that the approximate \ac{hac} \eqref{TSG-eq: approximation 4} coincides with the dual-port \ac{gfm} control \cite{gross2022dual}.
\subsection{LC filter consideration and cascaded controllers}
\begin{figure}[b!]
	\includegraphics[width=\columnwidth]{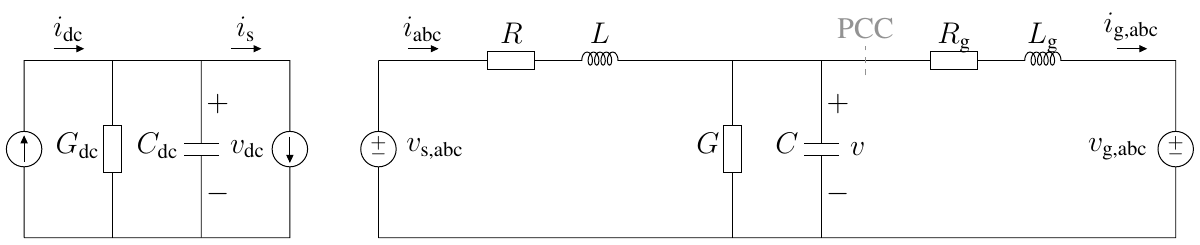}
	\caption{The circuit diagram of the converter model with LC output filter connected to a weak grid model in abc-coordinates system; this model also represents the case of a LCL filter consideration in which the grid-side filter inductance is merged with the grid impedance.}\label{TSG-fig:4}
\end{figure}

\begin{figure*}[t!]
	\centering	
	\includegraphics[width=\textwidth]{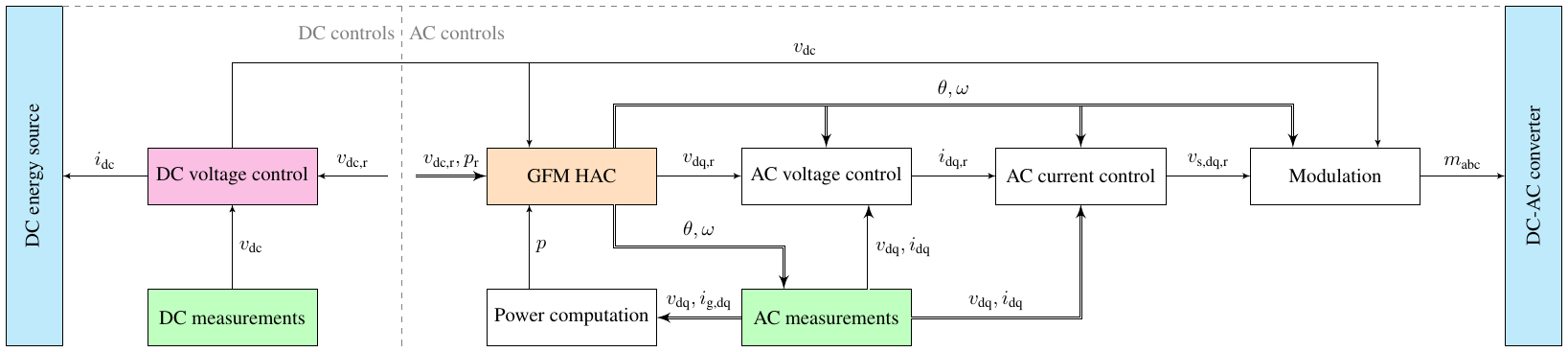}
	\caption{The overall control architecture of a grid-connected converter system including the dc-side controls, ac-side \ac{gfm}, voltage, and current controls.}\label{TSG-fig:5}
\end{figure*}

It is possible to consider LC output filter element which results in the model configuration in Figure \ref{TSG-fig:4}. In this case, one can combine \eqref{TSG-eq: approximation 3} with standard PI-based cascaded voltage and current controls \cite{TGAKD20,markovic2021understanding,subotic2020lyapunov}. Such control architecture is shown in Figure \ref{TSG-fig:5} and is briefly described as it follows.
\begin{itemize}
\item The phase angle defined by \ac{hac} and a prescribed reference ac voltage magnitude are combined to define the reference \ac{pcc} voltage in dq-coordinates, i.e., $v_\text{dq,r}$. Note that the converter frequency and angle defined by integrating \eqref{TSG-eq: approximation 3} serves as the reference angle for subsequent controllers implemented in dq0-coordinates.
 \item The PI-based ac voltage control (with feedforward terms) \cite{TGAKD20} processes the reference given by the \ac{gfm} layer and \ac{pcc} voltage feedback to define the reference filter current, i.e., $i_\tx{dq,r}$.
 \item The PI-based ac current control (with feedforward terms) processes the reference given by the voltage control layer and filter current feedback to define the converter voltage to appear behind the LC filter, i.e., $v_\tx{s,dq,r}$.
 \item The reference voltage given by the current control is processed by the modulation algorithm that defines $m_\tx{abc}$.   
\end{itemize}
We refer the reader to \cite{TGAKD20,subotic2020lyapunov,markovic2021understanding} for further details on the structure and tuning of such cascaded control architectures. 
\subsection{Control variants and extensions}
On the basis of the \ac{gfm} \ac{hac}, one can construct a few control variants. Let us introduce three different variants.  
\subsubsection{Fully multi-variable variant} the key idea behind the \ac{hac} is to include a dc feedback controller into the converter angle dynamics. Along the same direction, one can include an ac feedback controller into the converter dc dynamics. In a generic form, the closed-loop dynamics under such fully multi-variable control design takes the form 
	\begin{align*}
	\ddt{x_\tx{dc}}&=f_\tx{dc}(x_\tx{dc},x_\tx{ac})+\kappa_\tx{dc}g_{11}(x_\tx{dc})+\kappa_{\tx{ac}\to\tx{dc}}g_{12}(x_\tx{ac}),
	\\
	\ddt{x_\tx{ac}}&=f_\tx{ac}(x_\tx{dc},x_\tx{ac})+\kappa_{\tx{dc}\to\tx{ac}}g_{21}(x_\tx{dc})+\kappa_\tx{ac}g_{22}(x_\tx{ac}),	
	\end{align*}
	where $x_\tx{dc}$ denotes the dc states, $f_\tx{dc}(x_\tx{dc},x_\tx{ac})$ describes the physical dc subsystem, $\kappa_\tx{dc}$ is the dc control gain, $g_{11}(x_\tx{dc})$ is the linear/nonlinear dc controller for the dc states, $\kappa_{\tx{ac}\to\tx{dc}}$ is the gain of ac$\to$dc linear/nonlinear coupling control  $g_{12}(x_\tx{ac})$. Similarly, the states, physical ac subsystem, coupling controller gain and function, the ac control gain and function of the ac subsystem are respectively denoted by $x_\tx{ac}$, $f_\tx{ac}(x_\tx{dc},x_\tx{ac})$, $\kappa_{\tx{dc}\to\tx{ac}}$, $g_{21}(x_\tx{dc})$, $\kappa_\tx{ac}$, and $g_{22}(x_\tx{ac})$. Such augmentation of the standard \ac{gfm} \ac{hac} is explored in \cite{Andrea}.
\subsubsection{Inverse tangent variant} one can replace the ac synchronization term in \eqref{TSG-eq:6} by an inverse tangent function, i.e., $\tan^{-1}(\delta-\delta_\tx{r})$. Furthermore, one can consider the combination of controls in \cite{rowe2012arctan} with \cite{arghir2019electronic} to arrive at the hybrid form:
	\begin{equation*}
	\omega=\omega_0+\kappa_\tx{dc}\left(v_\tx{dc}-v_\tx{dc,r}\right)-\kappa_\tx{ac,1}\tan^{-1}\left(\kappa_\tx{ac,2}(p-p_\tx{r})\right).
	\end{equation*}
Note that it is possible to derive strong large-signal (i.e., global) stability guarantees for this control variant as in \cite{tayyebi2022grid,Andrea}. Further, \cite{rowe2012arctan} highlights the improved dynamic performance of the arctan droop control in contrast to the standard droop control. Thus, one can expect similar improvements for the hybrid arctan variant versus the standard \ac{hac}.
\subsubsection{Energy-like variant} finally, one can replace the linear dc term in \eqref{TSG-eq:6} with a nonlinear quadratic term, i.e., $\left(v_\tx{dc}-v_\tx{dc,r}\right)^2$ that is related to the dc energy and its reference. This control variant is particularly interesting for the modular multilevel converter (MMC) applications. 

We remark that our preliminary investigations, e.g., \cite{Andrea}, suggest that the aforementioned variants exhibit improved performance and/or lead to more relaxed conditions over the standard \ac{hac}, however, a deeper investigation is required.  
\section{Experimental verification}\label{TSG-sec:verification}
In this section, we describe the employed \ac{chil} verification approach, and present our test results.
\subsection{Controller-hardware-in-the-loop verification approach}
To verify the proposed \ac{gfm} \ac{hac} strategy under real-world effects such as discretization, delays, measurement noise, etc., we go beyond offline simulations as in \cite{tayyebi2020almost,tayyebi2022grid,tayyebi2022system,tayyebi2022hybrid} and implement our control algorithm in a control card, in order to run \ac{chil} simulations. The \ac{chil} approach represents a good candidate in terms of balancing testing complexity, costs, and fidelity. This setup enables a high degree of automation, thereby facilitating a high coverage of cases and grid conditions, especially those hard to implement in a laboratory setup or in the field.

\begin{figure}[b!]
    \centering
	\includegraphics[width=\columnwidth]{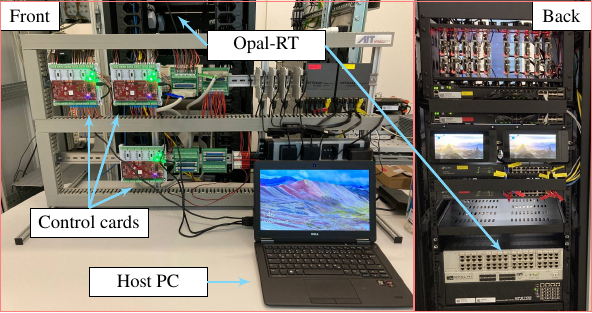}
	\caption{\ac{chil} testbed that includes three external control cards connected to Opal-RT OP5700. The host PC is used to automatically run \ac{chil} simulations and collect the results.}\label{TSG-fig:hardware}
    \vspace{3mm}
    \includegraphics[width=\columnwidth]{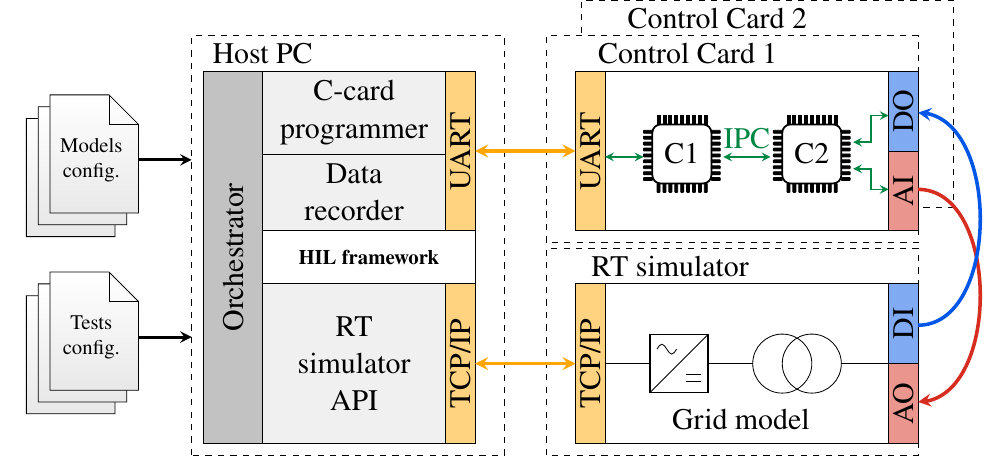}
    \caption{\ac{chil} testbed architecture. The \ac{chil} framework (on the left side) is developed in Python and its main functionality is to coordinate the RT-simulator and control cards. The control cards (on the right side) are connected to the RT-simulator via physical cables.}\label{TSG-fig:hil_architecture}
\end{figure}

The utilized hardware benchmark  is depicted in Figure \ref{TSG-fig:hardware}, consisting of an Opal-RT OP5700 as real-time simulator, a host PC, and several Texas Instrument (LaunchPad F28379D) control cards  in charge of executing the controller.  The control cards receive the dc voltage and ac voltage and current as analog signals from Opal-RT, and generate the PWM signals to be sent back to the real-time simulator as digital signals. Switching frequency for the inverter is set to 5kHz, which is the same rate as for the execution of the controller in the control card. The grid-connected converter model and the dc source controller are executed in Opal-RT, using Time-Stamped Bridges to model the inverter IGBTs \cite{dufour2005real}.

The testbed architecture is represented in Figure \ref{TSG-fig:hil_architecture}. The host PC communicates with the control cards via UART, and with Opal-RT via TCP/IP, thanks to the RT-LAB API. Given that the control cards possess two cores, data recording occurs in an online manner, bypassing memory limitations in the control card. Our setup allows us to measure relevant internal signals from the controllers such as frequency and filtered ac power. By means of configuration files, the tests of interest, set-points, and models to be used are defined. The tests are completely automatized, including the flashing of the control cards, building the grid models for the Opal-RT, synchronization of the cards and the real-time simulator, and finally the retrieval of all data of interest. 
\subsection{\ac{chil} verification test scenarios}
In what follows, we provide the results of four verification test cases that are performed on the testbed shown in Figure~\ref{TSG-fig:hardware}. Let us begin by highlighting the combination of employed controllers. We consider the PI-based dc and ac voltage control as described in Section \ref{TSG-sec: stiff grid}, and the approximate power-based implementation of the \ac{hac} presented in Section \ref{TSG-sec: weak grid}. The baseline grid-connected converter model that is implemented in Opal-RT corresponds to the circuit configuration illustrated in Figure \ref{TSG-fig:4}. Finally, the baseline model and control parameters are presented in Table \ref{TSG-table:1}. Note that the test-specific model and parameters modifications are described case-by-case.
\begin{table}[b!]
	\begin{center}
		\caption{Grid-connected converter model and control parameters.}\label{TSG-table:1}
		\begin{tabular}{ c | c | c }
			\hline\hline
			Symbol & Description &   Value \\
			\hline\hline
			$p_\tx{b}$ & base power & 500 kVA \\ 
			$f_\tx{b}$ & base frequency & 60 Hz \\ 
			$\omega_0$ & reference angular frequency & 2$\pi f_\tx{b}$ \\ 
			$v_0$ & reference grid voltage magnitude & 326.59 V\\ 
			$v_\tx{dc,r}$ & reference dc-link voltage & 3$v_0$\\ 
			$G_\tx{dc}$ & dc-side conductance & 0.01 m$\Omega^{-1}$\\ 
			$C_\tx{dc}$ & dc-link capacitance & 0.01 F\\ 
			$f_\tx{sw}$ & switching frequency & 5 kHz\\ 
			$L$ & ac filter inductance & 0.12 mH\\ 
			$C$ & ac filter capacitance & 0.13 mF\\ 
			$L_\tx{g}$ & grid equivalent inductance & 0.56 mH\\ 
			$R_\tx{g}$ & grid equivalent resistance & 0.064 $\Omega$\\ 
			$\kappa_\tx{p}$ & dc voltage control proportional gain & 10\\
			$\kappa_\tx{i}$ & dc voltage control integral gain & 500\\
			$\kappa_\tx{p,ac}$ & PCC voltage control proportional gain & 0.1\\
			$\kappa_\tx{i,ac}$ & PCC voltage control integral gain & 20\\
			$\kappa_\tx{dc}$ & HAC dc gain & 0.18\\
			$\bar{\kappa}_\tx{ac}$ & HAC ac gain & 18.84\\
			\hline\hline
		\end{tabular}
	\end{center}
\end{table}
\subsubsection{Accuracy verification and islanded \ac{gfm} operation} in this test scenario, the grid model is removed from the configuration in Figure \ref{TSG-fig:4}. Instead a resistive load is connected at the \ac{pcc} which at rated ac voltage consumes $\mathrm{0.5~p.u.}$ active power. Figure \ref{TSG-fig:islanded} illustrates the behavior of the islanded converter under the \ac{hac} in offline and \ac{chil} simulations. Note that a $\mathrm{0.5~p.u.}$ load increase is applied at $\mathrm{t=0.1~s}$. The results of offline and \ac{chil} simulations are sufficiently close, thus, verifying the accuracy of \ac{chil} testbed. Furthermore, the dynamic behavior shown in Figure \ref{TSG-fig:islanded} verifies the performances of the \ac{hac} control in islanded configuration. Note that, the ac gain (i.e., the droop gain) is selected such that it results in $\mathrm{5\%}$ frequency deviation for $\mathrm{1~p.u.}$ active power disturbance. Observe that the $\mathrm{0.5~p.u.}$ results in $\mathrm{2.5\%}$ frequency drop in Figure \ref{TSG-fig:islanded}, thus, verifying the drooping behavior of the approximate \ac{hac} \eqref{TSG-eq: approximation 3}. Finally, the dc voltage is recovered to the reference value due to the integral term.   
\subsubsection{Grid-connected \ac{gfm} operation} in this scenario, model configuration is identical to the one showed in Figure \ref{TSG-fig:4} and \ac{hac} behavior is investigated \ac{wrt} a set-point change event in grid-connected mode. Figure \ref{TSG-fig:grid-connected} shows that \ac{gfm} \ac{hac} not only preserves synchronization with the grid under a relatively large active power set-point change, i.e., $\mathrm{0.5~p.u.}$ increase, but also achieves zero post-event stead-state error and stabilizes the converter frequency at the desired reference.  We remark that the difference in transient behaviors in Figures \ref{TSG-fig:islanded} and \ref{TSG-fig:grid-connected} originates from the natural damping influence of the resistive load in the previous test scenario. Finally, retuning the \ac{gfm} control parameters and the cut-off frequency of the low-pass filter applied to ac power measurement allows to realize a first-order behavior following the set-point change event. However, for the sake of consistency the parameters are kept identical. 
\subsubsection{Grid frequency variation} in this scenario, the grid frequency is increased step-wise by $\mathrm{5\%}$. Figure \ref{TSG-fig:frequency-change} highlights the behavior of grid-connected converter under the \ac{hac}. Note that the \ac{gfm} \ac{hac} preserves system stability \ac{wrt} a relatively severe contingency. This is underpinned by the large-signal (i.e., global) stability of the \ac{hac} established in \cite{tayyebi2022grid}. Further, due to the particular choice of the droop gain, the converter active power injection drops by $\mathrm{1~p.u.}$ and reverses the power flow to provide frequency support.
\subsubsection{Two-converter load-sharing operation} 
\begin{figure*}[t!]
	\centering	
	{\includegraphics[trim=8mm 6mm 16mm 6mm,clip,width=0.32\textwidth]{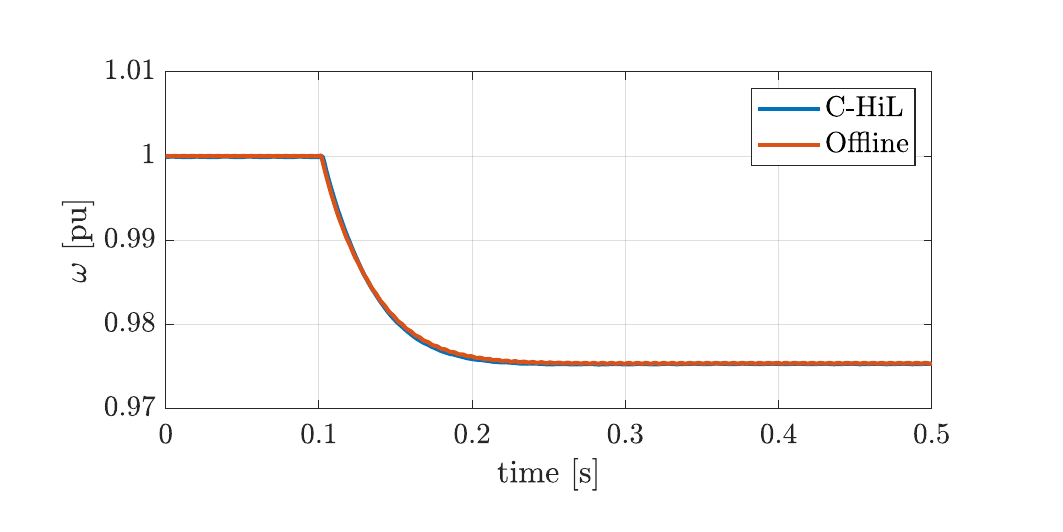}}
	\includegraphics[trim=8mm 6mm 16mm 6mm,clip,width=0.32\textwidth]{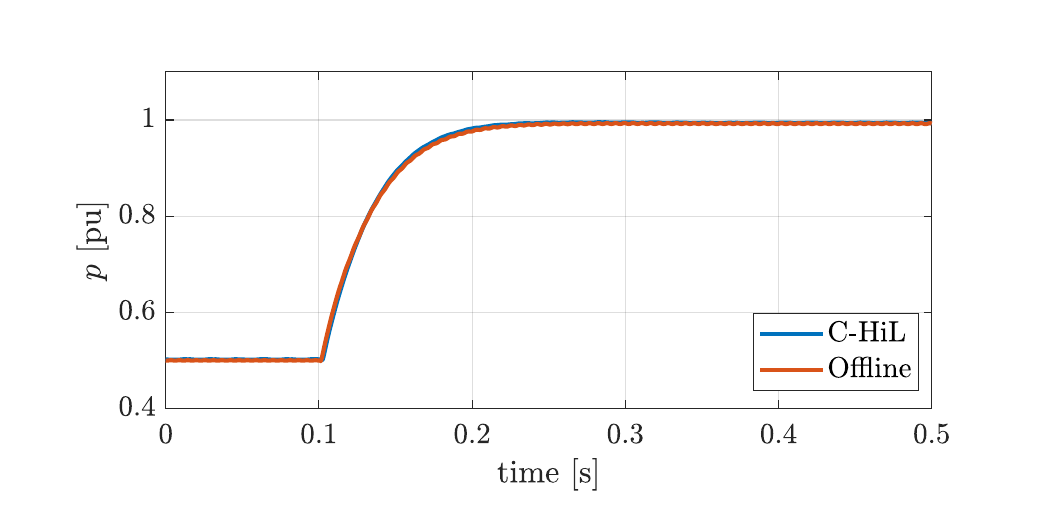}
	\includegraphics[trim=8mm 6mm 16mm 6mm,clip,width=0.32\textwidth]{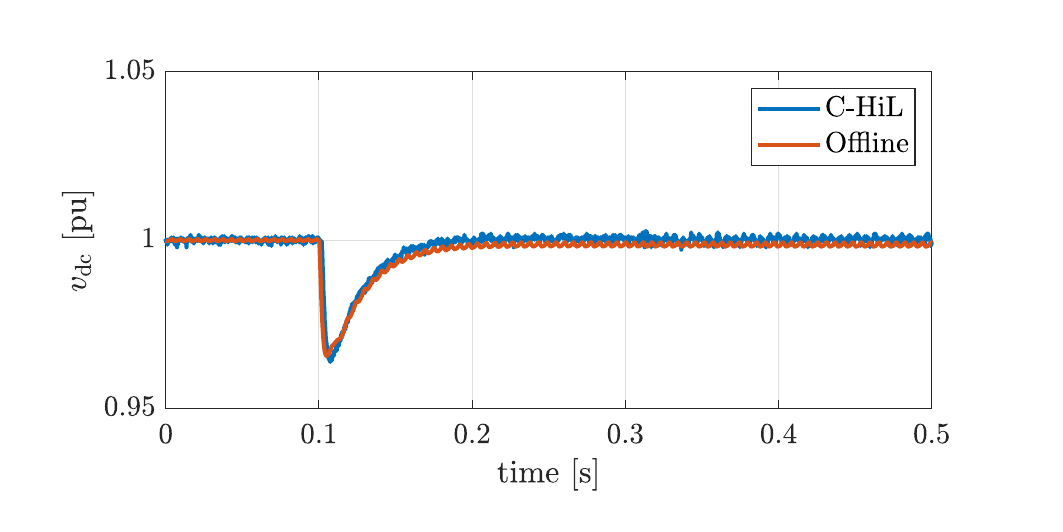}
	\caption{The time-evolution of normalized frequency (left), active power (middle), and dc voltage (right) of an islanded \ac{gfm} converter under \ac{hac} \ac{wrt} a load disturbance scenario in offline and \ac{chil} simulations.}\label{TSG-fig:islanded}\vspace{3mm}
	{\includegraphics[trim=8mm 6mm 16mm 6mm,clip,width=0.32\textwidth]{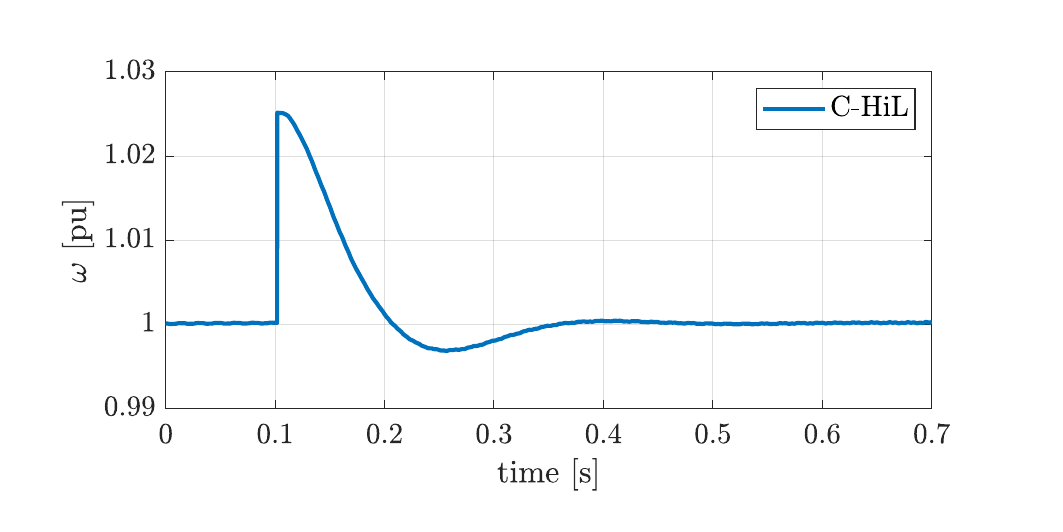}}
	\includegraphics[trim=8mm 6mm 16mm 6mm,clip,width=0.32\textwidth]{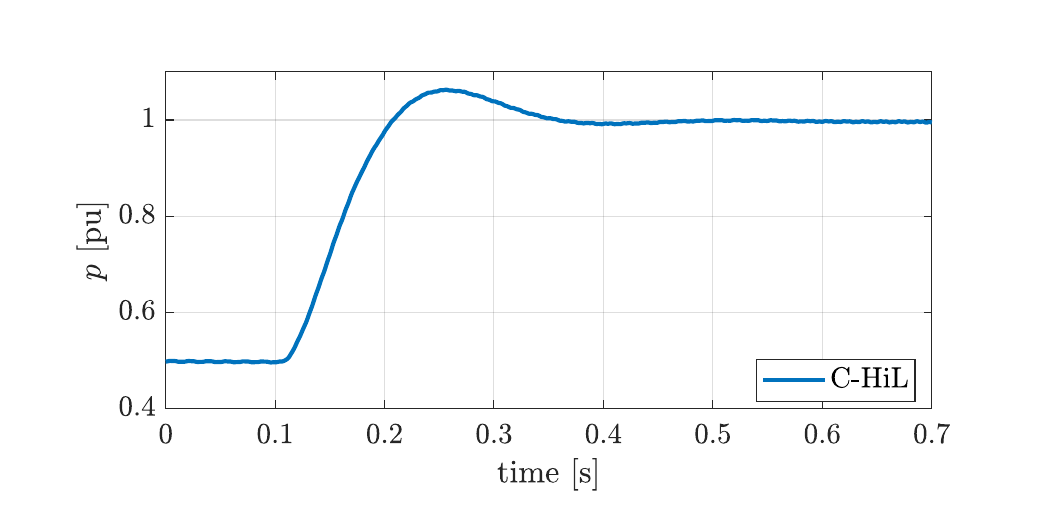}
	\includegraphics[trim=8mm 6mm 16mm 6mm,clip,width=0.32\textwidth]{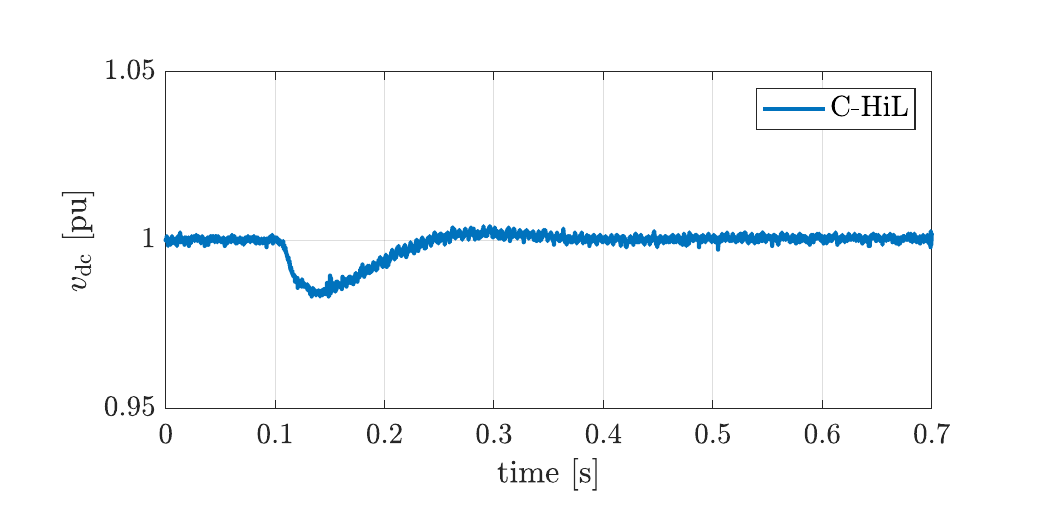}
	\caption{The time-evolution of normalized frequency (left), active power (middle), and dc voltage (right) of a grid-connected \ac{gfm} converter under \ac{hac} \ac{wrt} a power set-point change scenario in \ac{chil} simulations.}\label{TSG-fig:grid-connected}\vspace{3mm}
	{\includegraphics[trim=8mm 6mm 16mm 6mm,clip,width=0.32\textwidth]{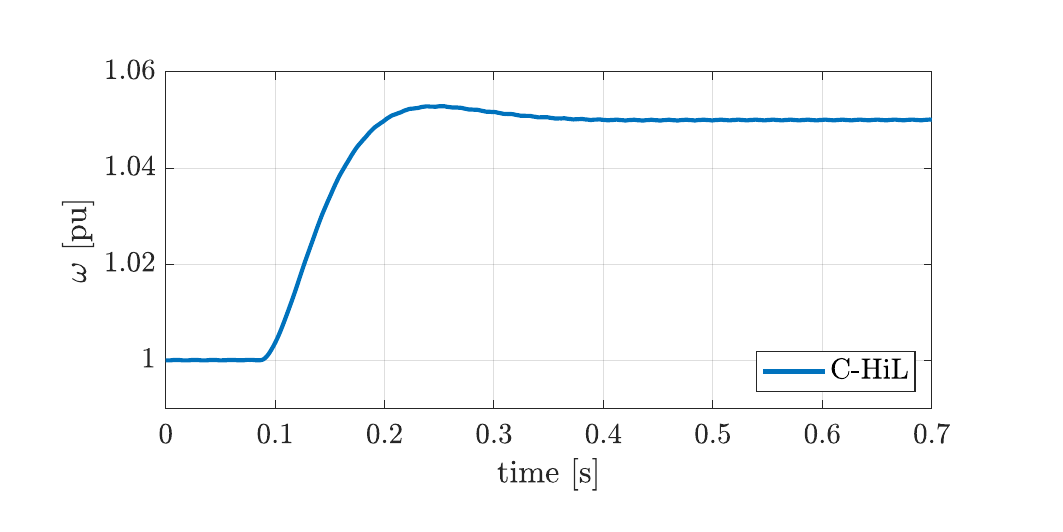}}
	\includegraphics[trim=8mm 6mm 16mm 6mm,clip,width=0.32\textwidth]{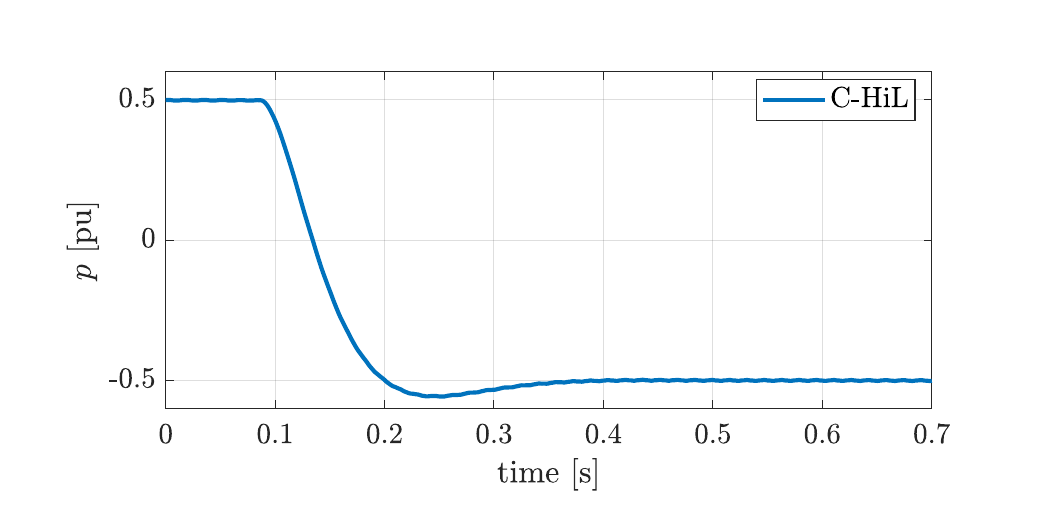}
	\includegraphics[trim=8mm 6mm 16mm 6mm,clip,width=0.32\textwidth]{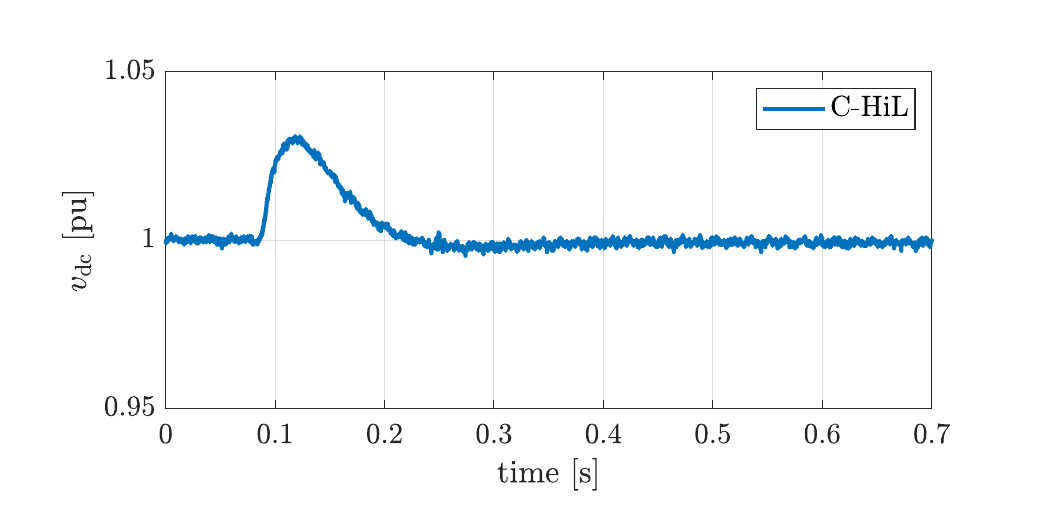}
	\caption{The time-evolution of normalized frequency (left), active power (middle), and dc voltage (right) of a grid-connected \ac{gfm} converter under \ac{hac} \ac{wrt} a grid frequency variation scenario in \ac{chil} simulations.}\label{TSG-fig:frequency-change}\vspace{3mm}
	{\includegraphics[trim=8mm 6mm 16mm 6mm,clip,width=0.32\textwidth]{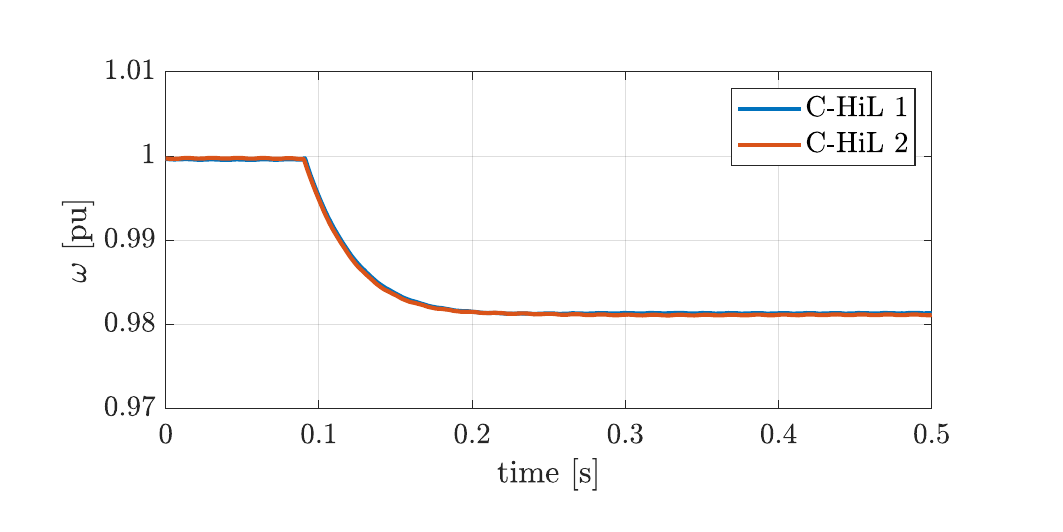}}
	\includegraphics[trim=8mm 6mm 16mm 6mm,clip,width=0.32\textwidth]{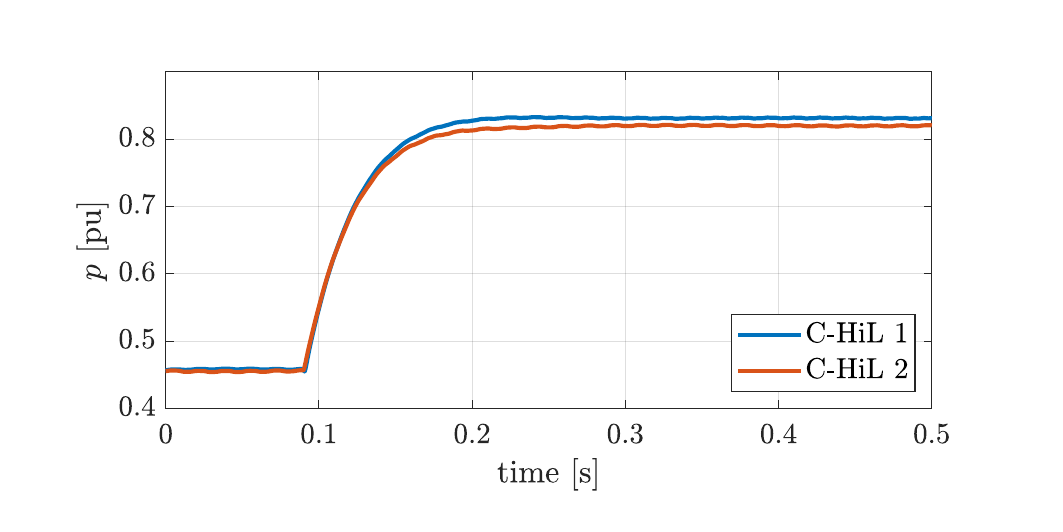}
	\includegraphics[trim=8mm 6mm 16mm 6mm,clip,width=0.32\textwidth]{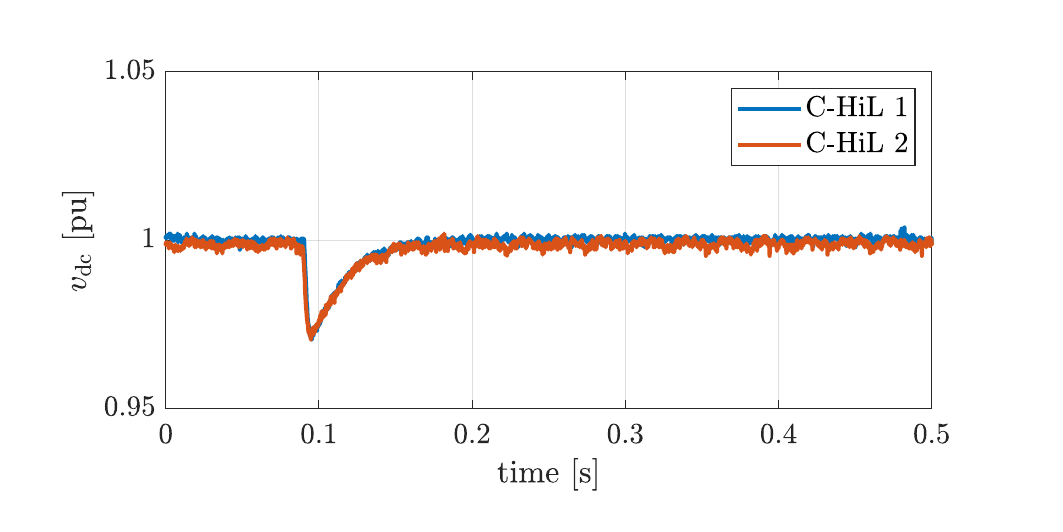}
	\caption{The time-evolution of normalized frequency (left), active power (middle), and dc voltage (right) of two coupled \ac{gfm} converters under \ac{hac} (with slightly different droop gains) \ac{wrt} a load disturbance scenario in \ac{chil} simulations.}\label{TSG-fig:two-converter}
\end{figure*}

finally, we consider a two-converter test scenarios. In this case, the converter models correspond to the model shown in Figure \ref{TSG-fig:4}. However, the grid model is removed and the converters are connected through two RL line models and a resistive load is connected in the middle. Note that the droop gains for the power converters are slightly different, i.e., $0.98\bar{\kappa}_\tx{ac}$ and $1.02\bar{\kappa}_\tx{ac}$. Figure \ref{TSG-fig:two-converter} illustrates the behavior of the system when a load increase event is applied. Observe that the post-disturbance frequency synchronization is achieved while the converters exhibits slightly different load-sharing according to the prescribed droop gains. 

\section{Conclusion}\label{TSG-sec:conclusion}
In this work, we discussed the behavioral properties of the \ac{gfm} \ac{hac}, described required complementary controls, provided a closed-loop analysis involving analytical operating point evaluation and energy-based nonlinear stability analysis, derived an approximation of the \ac{hac}, and introduced several extensions of the standard \ac{hac}. Last but not least, the control performance was verified by several \ac{chil} test scenarios. Our agenda of future work includes the stability analysis and performance verification of the \ac{hac} variants, and power hardware validation of the control concept.  
\bibliographystyle{IEEEtran}
\bibliography{IEEEabrv,Ref}

\begin{thebibliography}{10}
\providecommand{\url}[1]{#1}
\csname url@samestyle\endcsname
\providecommand{\newblock}{\relax}
\providecommand{\bibinfo}[2]{#2}
\providecommand{\BIBentrySTDinterwordspacing}{\spaceskip=0pt\relax}
\providecommand{\BIBentryALTinterwordstretchfactor}{4}
\providecommand{\BIBentryALTinterwordspacing}{\spaceskip=\fontdimen2\font plus
\BIBentryALTinterwordstretchfactor\fontdimen3\font minus
  \fontdimen4\font\relax}
\providecommand{\BIBforeignlanguage}[2]{{%
\expandafter\ifx\csname l@#1\endcsname\relax
\typeout{** WARNING: IEEEtran.bst: No hyphenation pattern has been}%
\typeout{** loaded for the language `#1'. Using the pattern for}%
\typeout{** the default language instead.}%
\else
\language=\csname l@#1\endcsname
\fi
#2}}
\providecommand{\BIBdecl}{\relax}
\BIBdecl

\bibitem{hatziargyriou2020definition}
N.~Hatziargyriou, J.~Milanovic, C.~Rahmann, V.~Ajjarapu, C.~Canizares,
  I.~Erlich, D.~Hill, I.~Hiskens, I.~Kamwa, B.~Pal \emph{et~al.}, ``Definition
  and classification of power system stability--revisited \& extended,''
  vol.~36, no.~4, pp. 3271--3281, 2020.

\bibitem{lin2020research}
Y.~Lin, J.~H. Eto, B.~B. Johnson, J.~D. Flicker, R.~H. Lasseter, H.~N.
  Villegas~Pico, G.-S. Seo, B.~J. Pierre, and A.~Ellis, ``Research roadmap on
  grid-forming inverters,'' National Renewable Energy Lab. (NREL), Tech. Rep.,
  2020.

\bibitem{crivellaro2020beyond}
A.~Crivellaro, A.~Tayyebi, C.~Gavriluta, D.~Gro{\ss}, A.~Anta, F.~Kupzog, and
  F.~D{\"o}rfler, ``Beyond low-inertia systems: Massive integration of
  grid-forming power converters in transmission grids,'' in \emph{{IEEE} Power
  \& Energy Society General Meeting (PESGM)}, 2020, pp. 1--5.

\bibitem{milano2018foundations}
F.~Milano, F.~D{\"o}rfler, G.~Hug, D.~J. Hill, and G.~Verbi{\v{c}},
  ``Foundations and challenges of low-inertia systems,'' in \emph{power systems
  computation conference (PSCC)}, 2018.

\bibitem{TGAKD20}
A.~{Tayyebi}, D.~{Groß}, A.~{Anta}, F.~{Kupzog}, and F.~{Dörfler},
  ``Frequency stability of synchronous machines and grid-forming power
  converters,'' \emph{{IEEE} Trans. Emerg. Sel. Topics Power Electron.},
  vol.~8, no.~2, pp. 1004--1018, 2020.

\bibitem{Rocabert}
J.~Rocabert, A.~Luna, F.~Blaabjerg, and P.~Rodríguez, ``Control of power
  converters in ac microgrids,'' \emph{{IEEE} Trans. Power Electron.}, vol.~27,
  no.~11, pp. 4734--4749, 2012.

\bibitem{markovic2021understanding}
U.~Markovic, O.~Stanojev, P.~Aristidou, E.~Vrettos, D.~Callaway, and G.~Hug,
  ``Understanding small-signal stability of low-inertia systems,'' vol.~36,
  no.~5, pp. 3997--4017, 2021.

\bibitem{NGESO2023}
``Great {B}ritain grid forming best practice guide,'' National Grid ESO, Tech.
  Rep., 2023.

\bibitem{AEMO2022}
``Engineering roadmap to 100\% renewables,'' Australian Energy Market Operator
  (AEMO), Tech. Rep., 2022.

\bibitem{CDA93}
M.~Chandorkar, D.~Divan, and R.~Adapa, ``Control of parallel connected
  inverters in standalone {AC} supply systems,'' \emph{{IEEE} Trans. Ind.
  Appl.}, vol.~29, no.~1, pp. 136--143, 1993.

\bibitem{rowe2012arctan}
C.~N. Rowe, T.~J. Summers, R.~E. Betz, D.~J. Cornforth, and T.~G. Moore,
  ``Arctan power--frequency droop for improved microgrid stability,''
  \emph{{IEEE} Trans. Power Electron.}, vol.~28, no.~8, pp. 3747--3759, 2012.

\bibitem{yu2020comparative}
H.~Yu, M.~Awal, H.~Tu, I.~Husain, and S.~Lukic, ``Comparative transient
  stability assessment of droop and dispatchable virtual oscillator controlled
  grid-connected inverters,'' \emph{{IEEE} Trans. Power Electron.}, vol.~36,
  no.~2, pp. 2119--2130, 2021.

\bibitem{ZW11}
Q.~C. Zhong and G.~Weiss, ``Synchronverters: inverters that mimic synchronous
  generators,'' \emph{{IEEE} Trans. Ind. Electron.}, vol.~58, no.~4, pp.
  1259--1267, 2011.

\bibitem{chen2021enhanced}
M.~Chen, D.~Zhou, and F.~Blaabjerg, ``Enhanced transient angle stability
  control of grid-forming converter based on virtual synchronous generator,''
  \emph{{IEEE} Trans. Ind. Electron.}, vol.~69, no.~9, pp. 9133--9144, 2021.

\bibitem{arghir2019electronic}
C.~Arghir and F.~D{\"o}rfler, ``The electronic realization of synchronous
  machines: Model matching, angle tracking, and energy shaping techniques,''
  \emph{{IEEE} Trans. Power Electron.}, vol.~35, no.~4, pp. 4398--4410, 2019.

\bibitem{cvetkovic_modeling_2015}
I.~Cvetkovic, D.~Boroyevich, R.~Burgos, C.~Li, and P.~Mattavelli, ``Modeling
  and control of grid-connected voltage-source converters emulating isotropic
  and anisotropic synchronous machines,'' in \emph{IEEE Workshop on Control and
  Modeling for Power Electronics (COMPEL)}, 2015.

\bibitem{aghdam2022virtual}
S.~A. Aghdam and M.~Agamy, ``Virtual oscillator-based methods for grid-forming
  inverter control: A review,'' \emph{IET Renew. Power Gener.}, vol.~16, no.~5,
  pp. 835--855, 2022.

\bibitem{seo2019dispatchable}
G.-S. Seo, M.~Colombino, I.~Subotic, B.~Johnson, D.~Gro{\ss}, and
  F.~D{\"o}rfler, ``Dispatchable virtual oscillator control for decentralized
  inverter-dominated power systems: Analysis and experiments,'' in \emph{IEEE
  Applied Power Electronics Conference and Exposition (APEC)}, 2019.

\bibitem{awal2022double}
M.~Awal, M.~R.~K. Rachi, H.~Yu, I.~Husain, and S.~Lukic, ``Double synchronous
  unified virtual oscillator control for asymmetrical fault ride-through in
  grid-forming voltage source converters,'' \emph{{IEEE} Trans. Power
  Electron.}, 2022.

\bibitem{tayyebi2022grid}
A.~Tayyebi, A.~Anta, and F.~D{\"o}rfler, ``Grid-forming hybrid angle control
  and almost global stability of the dc-ac power converter,'' \emph{IEEE
  Transactions on Automatic Control}, 2022.

\bibitem{gao2020grid}
Y.~Gao, H.-P. Ren, and J.~Li, ``Grid-forming converters control based on dc
  voltage feedback,'' 2020, {Preprint} available at
  {https://arxiv.org/abs/2009.05759}.

\bibitem{gross2022dual}
D.~Gro{\ss}, E.~S{\'a}nchez-S{\'a}nchez, E.~Prieto-Araujo, and
  O.~Gomis-Bellmunt, ``Dual-port grid-forming control of mmcs and its
  applications to grids of grids,'' \emph{IEEE Transactions on Power Delivery},
  vol.~37, no.~6, pp. 4721--4735, 2022.

\bibitem{chen2022generalized}
M.~Chen, D.~Zhou, A.~Tayyebi, E.~Prieto-Araujo, F.~D{\"o}rfler, and
  F.~Blaabjerg, ``Generalized multivariable grid-forming control design for
  power converters,'' \emph{{IEEE} Trans. Smart Grid}, vol.~13, no.~4, pp.
  2873--2885, 2022.

\bibitem{simpson2013synchronization}
J.~W. Simpson-Porco, F.~D{\"o}rfler, and F.~Bullo, ``Synchronization and power
  sharing for droop-controlled inverters in islanded microgrids,''
  \emph{Automatica}, vol.~49, no.~9, pp. 2603--2611, 2013.

\bibitem{tayyebi2022hybrid}
A.~Tayyebi and F.~D{\"o}rfler, ``Hybrid angle control and almost global
  stability of non-synchronous hybrid ac/dc power grids,'' in \emph{IEEE
  Conference on Decision and Control (CDC)}, 2022.

\bibitem{tayyebi2020almost}
A.~Tayyebi, A.~Anta, and F.~D{\"o}rfler, ``Almost globally stable grid-forming
  hybrid angle control,'' in \emph{IEEE Conference on Decision and Control
  (CDC)}, 2020.

\bibitem{tayyebi2022system}
A.~Tayyebi, A.~Magdaleno, D.~Vettoretti, M.~Chen, E.~Prieto-Araujo, A.~Anta,
  and F.~D{\"o}rfler, ``System-level performance and robustness of the
  grid-forming hybrid angle control,'' \emph{Electric Power Systems Research},
  vol. 212, p. 108503, 2022.

\bibitem{arghir2018grid}
C.~Arghir, T.~Jouini, and F.~D{\"o}rfler, ``Grid-forming control for power
  converters based on matching of synchronous machines,'' \emph{Automatica},
  vol.~95, pp. 273--282, 2018.

\bibitem{yazdani10voltage}
A.~Yazdani and R.~Iravani, \emph{Voltage-sourced converters in power systems:
  modeling, control, and applications}.\hskip 1em plus 0.5em minus 0.4em\relax
  John Wiley \& Sons, 2010.

\bibitem{colon2022stability}
G.~E. Col{\'o}n-Reyes, K.~C. Stocking, D.~S. Callaway, and C.~J. Tomlin,
  ``Stability and robustness of a hybrid control law for the half-bridge
  inverter,'' 2022, preprint available at: https://arxiv.org/abs/2204.07539.

\bibitem{albea2017hybrid}
C.~Albea, O.~L. Santos, D.~Z. Prada, F.~Gordillo, and G.~Garcia, ``Hybrid
  control scheme for a half-bridge inverter,'' \emph{IFAC-PapersOnLine},
  vol.~50, no.~1, pp. 9336--9341, 2017.

\bibitem{CGD17}
S.~Curi, D.~Gro\ss, and F.~D\"{o}rfler, ``Control of low-inertia power grids: a
  model reduction approach,'' in \emph{{IEEE} {Conference} on {Decision} and
  {Control} ({CDC})}, 2017.

\bibitem{samanta2022fast}
S.~Samanta, N.~R. Chaudhuri, and C.~Lagoa, ``Fast frequency support from
  grid-forming converters under dc-and ac-side current limits,'' \emph{IEEE
  Transactions on Power Systems}, 2022.

\bibitem{samanta2021stability}
S.~Samanta and N.~R. Chaudhuri, ``Stability analysis of grid-forming converters
  under dc-side current limitation in primary frequency response regime,''
  \emph{IEEE Transactions on Power Systems}, vol.~37, no.~4, pp. 3077--3091,
  2021.

\bibitem{ulbig2014impact}
A.~Ulbig, T.~S. Borsche, and G.~Andersson, ``Impact of low rotational inertia
  on power system stability and operation,'' in \emph{IFAC Proceedings
  Volumes}, 2014.

\bibitem{Andrea}
A.~Gattiglio, ``Multi-variable arctan hybrid angle control and global stability
  of grid-forming power converters,'' 2021, {Thesis} available at:
  http://dx.doi.org/10.13140/RG.2.2.18935.75684.

\bibitem{subotic2020lyapunov}
I.~Suboti{\'c}, D.~Gro{\ss}, M.~Colombino, and F.~D{\"o}rfler, ``A lyapunov
  framework for nested dynamical systems on multiple time scales with
  application to converter-based power systems,'' \emph{{IEEE} Trans. Autom.
  Control}, vol.~66, no.~12, pp. 5909--5924, 2020.

\bibitem{dufour2005real}
C.~Dufour and J.~B{\'e}langer, ``Real-time simulation of a 48-pulse {GTO
  STATCOM} compensated power system on a {Dual-Xeon PC} using {RTLAB},'' in
  \emph{International Conference on Power Systems Transients}, 2005.

\end{thebibliography}
\end{document}